\documentclass[number,citesort,dvips]{arxbj}
\usepackage{dcolumn}
\usepackage{graphicx}

\aid{0}
\volume{17}
\issue{2}
\pubyear{2011}
\firstpage{781}
\lastpage{813}
\doi{10.3150/10-BEJ293}

\makeatletter
\newcommand{\eqref}[1]{(\ref{#1})}
\newproclaim{assumption}{Assumption}
\newproclaim{Definition}[Theorem]{Definition}

\newtheorem{Lemma}[Theorem]{Lemma}
\newtheorem{Proposition}[Theorem]{Proposition}
\newtheorem{Theorem}{Theorem}[section]

\newremark{Remark}[Theorem]{Remark}
\newremark{example}{Example}[section]
\newcolumntype{d}[1]{D{.}{.}{#1}}
\makeatother

\begin{document}
\begin{frontmatter}

\title{Nonparametric tests for pathwise properties of semimartingales}
\runtitle{Nonparametric tests for pathwise properties of semimartingales}

\begin{aug}
\author[a]{\fnms{Rama} \snm{Cont}\corref{}\ead[label=e1]{Rama.Cont@columbia.edu}\thanksref{a}}
\and
\author[b]{\fnms{Cecilia} \snm{Mancini}\thanksref{b}\ead[label=e2]{cecilia.mancini@dmd.unifi.it}}

\runauthor{R. Cont and C. Mancini}

\address[a]{IEOR Department,
Columbia University,
New York and Laboratoire de Probabilit\'es et Mod\`eles Al\'eatoires,
CNRS-Universit\'e Paris VI,
France.
\printead{e1}}

\address[b]{Dipartimento di Matematica per le Decisioni,
Universit\`{a} di Firenze,
Italy.\\
\printead{e2}}
\end{aug}

\received{\smonth{5} \syear{2009}}
\revised{\smonth{3} \syear{2010}}

%
\begin{abstract}
We propose two nonparametric tests for
investigating the pathwise properties of a signal
modeled as the sum of a L\'evy process and a
Brownian semimartingale. Using a nonparametric threshold estimator
for the continuous component of the quadratic variation, we design
a test for the presence of a continuous
martingale component in the process and a test for establishing
whether the jumps have finite or infinite variation, based on
observations on a discrete-time grid. We evaluate the performance
of our tests using simulations of various stochastic models and use
the tests to investigate the fine structure of the DM/USD exchange
rate fluctuations and SPX futures prices. In both cases, our tests
reveal the presence of a non-zero Brownian component and a finite
variation jump component.
\end{abstract}

%

\end{frontmatter}

\section{Introduction}
Continuous-time stochastic models based on \textit{discontinuous
semimartingales} have been increasingly used in many applications,
such as financial econometrics, option pricing and stochastic
control.
Some of these models are constructed by
adding i.i.d. jumps to a continuous process driven by Brownian motion
\cite{merton76,kou}, while others are based on purely discontinuous
processes which move only through jumps \cite{Madan01,cgmy02}. Even
within the class of purely discontinuous models, one finds a variety
of models with different path properties -- finite/infinite jump
intensity, finite/infinite variation -- which turn out to have an
importance in applications, such as optimal stopping \cite{alili04}
and
the asymptotic behavior of option prices \cite{carrwu03,contankov}.
It is therefore of interest to investigate which class of models
-- diffusion, jump-diffusion or pure-jump -- is the most appropriate
for a given data set. Nonparametric procedures have been recently
proposed for investigating the presence of jumps
\cite{bnsBPV06a,ASJ06,leemykland07} and studying some fine
properties of the jumps
\cite{ASJ07,AitJacod08FavsIA,TodTau08,WoernerHatAlpha06b} in a
signal.
Here, we
address related, but different, issues: for a semimartingale
whose jump component is a L\'evy process, we propose a test for the
presence of a continuous martingale component in the price process,
which allows us to discriminate between pure-jump and jump-diffusion
models, and a test for determining whether the jump component has
finite or infinite variation. Our tests are based on a
nonparametric threshold estimator \cite{mancini07} for the
integrated variance (defined as the continuous component of the
quadratic variation) based on observations on a discrete-time grid.
Without imposing restrictive assumptions on the continuous
martingale component, we obtain a central limit theorem for this
threshold estimator (Section \ref{CLT.sec}) and use it to design our
tests (Section \ref{testDescription.sec}).

Using simulations of stochastic models commonly used in finance, we
check the performance of our tests for realistic sample sizes
(Section \ref{simuls.sec}). Applied to time series of the DM/USD
exchange rate and SPX futures prices (Section \ref{data.sec}), our
tests reveal, in both cases, the presence of a non-zero Brownian
component, combined with a finite variation jump component. These
results suggest that these asset prices may be modeled as the sum of
a Brownian martingale and a jump component of finite variation.

\section{Definitions and notation}\label{defs.sec}
We consider a semimartingale $(X_{t})_{t\in[0,T]}$, defined on a
(filtered) probability space $(\Omega, (\mathcal{F}_t)_{t\in[0, T]},
\mathcal{F},$ $P)$ with paths in $D([0,T],\mathbb{R})$, driven by a (standard)
Brownian motion $W$ and a pure-jump L\'evy process $L$:
\begin{equation} \label{eqX}
X_t= x_0+\int_0^t a_s\,\mathrm{d}s+\int_0^t\sigma_s \,\mathrm{d}W_s+ L_t,\qquad
  t\in\,]0,T],
\end{equation}
where $a$, $\sigma$ are adapted processes with
right-continuous paths with left limits (cadlag processes), such that
(\ref{eqX}) admits a unique strong solution $X$
on $[0,T]$ which is adapted and cadlag \cite{IkeWat81}. $L$~has L\'evy
measure $\nu$ and may be decomposed as
$L_t =
J_t + M_{t}$, where
\begin{equation} \label{RepresJu} J_{t}:=\int_0^t\!\!\!\int_{|x|> 1} x \mu
(\mathrm{d}x, \mathrm{d}s)=\sum_{\ell=1}^{N_t} \gamma_{\ell},\qquad   M_{t}:=
\int_0^t\!\!\!\int_{|x|\leq1} x[\mu(\mathrm{d}x, \mathrm{d}s)-\nu(\mathrm{d}x) \,\mathrm{d}t].
\end{equation}
$J$ is a
compound Poisson process representing the ``large'' jumps of $X$,
$\mu$ is a Poisson random measure on $[0,T]\times\mathbb{R}$ with
intensity measure $ \nu(\mathrm{d}x)\,\mathrm{d}t$, $N$ is a Poisson process with
intensity $\nu(\{x, |x|> 1\})<\infty$,
$\gamma_{\ell}$ are i.i.d. and independent of $N$ and the martingale $M$
is the compensated sum of small jumps of $L$. We will define $ \mu(\mathrm{d}x,
\mathrm{d}t)-\nu(\mathrm{d}x) \,\mathrm{d}t=:\tilde\mu( \mathrm{d}x, \mathrm{d}t)$, the compensated Poisson random
measure associated to $\mu$.
We allow for the \textit{infinite activity} (IA) case $\nu(\mathbb
{R})=\infty$, where small
jumps of $L$ occur infinitely often.
For a~semimartingale $Z,$ we denote by
$\Delta_{i} Z=Z_{t_i}-Z_{t_{i-1}}$ its increments and by
$\Delta Z_t=Z_t - Z_{t-}$ its jump at time $t$.
The \textit{Blumenthal--Getoor (BG) index} of $L$, defined as
\[
\alpha:= \inf\biggl\{\delta\geq0,  \int_{|x|\leq1} |x|^\delta\nu
(\mathrm{d}x)< +
\infty\biggr\}\leq2,
\]
measures the degree of \textit{activity} of small jumps. A
compound Poisson process has $\alpha=0$, while an
$\alpha$-stable process has BG index equal to
$\alpha\in\,]0,2[$. The gamma process and the variance gamma (VG)
process are examples of infinite activity L\'evy processes with
$\alpha=0$. A pure-jump L\'{e}vy process with BG index $\alpha<1$ has
paths with \textit{finite variation}, while for $\alpha> 1$, the
sample paths have \textit{infinite variation} a.s. When
$\alpha=1$, the paths may have either finite or infinite variation
\cite{BSW06}. The normal inverse Gaussian process (NIG) and the
generalized hyperbolic L\'evy motion (GHL) have infinite variation
and $\alpha=1$. Tempered stable processes \cite{cgmy02,contankov}
allow for $\alpha\in[0,2[$. We call $\mathit{IV}=\int_0^T \sigma^2_u \,\mathrm{d}u$
the \textit{integrated variance} of $X$ and $\mathit{IQ}=\int_0^T \sigma^4_u \,\mathrm{d}u$
the \textit{integrated quarticity} of $X$, and we write
\[
X_{0 t}=\int_0^t a_s \,\mathrm{d}s + \int_0^t \sigma_s \,\mathrm{d}W_s,\qquad  X_{1 t}= X_{0
t}+ J_t.
\]
We will use the following assumption. 
%
%
\renewcommand{\theassumption}{\Alph{assumption}\arabic{assumption}}
\setcounter{assumption}{0}
\begin{assumption}\label{asA1}
\begin{equation}\label{A1}
\exists\alpha\in[0,2]\qquad \int_{|x|\leq\varepsilon} x^2 \nu
(\mathrm{d}x) \sim
\varepsilon^{2-\alpha}\qquad \mbox{as $\varepsilon\rightarrow0$,}
\end{equation}
where
$f(h)\sim g(h)$ means that $f(h)=\mathrm{O}(g(h) )$ and $g(h)=\mathrm{O}(f(h) )$ as
$h\to0$.
\end{assumption}

This assumption implies that $\alpha$ is the BG index of
$L$. \textup{\ref{asA1}} is satisfied if, for instance,
$\nu$ has a density which behaves as $\frac{K_\pm}{|x|^{1+\alpha}}$
when $x\rightarrow0\pm$, where $K_\pm>0$. In particular, \textup{\ref{asA1}} holds
for all L\'evy processes commonly used in finance \cite{contankov}:
NIG, variance gamma,
tempered stable processes or generalized hyperbolic processes.

Typically, we observe $X_t$ in the form of a discrete record
$\{x_{0}, X_{t_1},\ldots,$ $X_{t_{n-1}}, X_{t_n}\}$
on a time grid
$t_i=ih$ with $h=T/n$. Our goal is to provide, given such a discrete
observations, nonparametric tests
for:
\begin{itemize}
\item detecting the presence of a continuous martingale
component in the price process;
\item analyzing the qualitative nature of the jump component, that is,
whether it has finite or infinite
variation.
\end{itemize}

%
\section{Central limit theorem for a threshold estimator of integrated
variance}\label{CLT.sec}
The ``realized variance'' $\sum_{  i=1}^n  ( \Delta_i X)^{ 2} $
of the semimartingale $X$ converges in probability \cite{Protter05}
to
\[
[X]_T:= \int_0^T \sigma_t^2 \,\mathrm{d}t + \int_0^T\!\!\!\int_{\mathbb{R}-\{0\}} x^2
\mu(\mathrm{d}x, \mathrm{d}s).
\]
A \textit{threshold estimator} \cite{mancini04,mancini07} of the
integrated variance $\mathit{IV}=\int_0^T \sigma_t^2 \,\mathrm{d}t$ is based on the idea of
summing only some of the squared increments of $X$, those
whose absolute value is smaller than some \textit{threshold} $r_h$:
\begin{equation} \label{hatIV} \hat{\mathit{IV}}_h := \sum_{i=1}^n (\Delta
_{i} X)^2I_{\{ (\Delta_{i} X)^2\leq r_h\}} .
\end{equation}
The term
$\int_0^T\!\!\!\int_{\mathbb{R}-\{0\}} x^2 \mu(\mathrm{d}x, \mathrm{d}s)$, due to jumps,
vanishes as
$h\to0$ for an appropriate choice of the threshold. P. L\'evy's
law for the modulus of continuity of the Brownian paths implies
that
\[
P\biggl(\lim _{h\rightarrow0}\sup_{i\in\{ 1,\ldots,n \}
}\frac{|\Delta_{i} W
|}{\sqrt{2h\ln 1/h}}
\leq1\biggr)=1
\]
and allows such a threshold to be chosen. It is shown in \cite
{mancini07}, Corollary 2, Theorem 4, that, under the above assumptions,
if we choose a
deterministic threshold $r_h$ such that
\begin{equation}\lim _{h\rightarrow0}r_h= 0 \quad \mbox{and}\quad
\lim _{h\rightarrow0}\frac{h\ln h}{ r_h}=
0,\label{rh.eq}
\end{equation}
then $\hat{\mathit{IV}}_h \stackrel{P}{\rightarrow} \mathit{IV}
$ as $h\rightarrow0$. If the jumps have finite intensity, then the thresholding
procedure allows as $h\to0$, a jump to be detected in $]t_{i-1},
t_i]$. In fact, since $a$ and $\sigma$ are cadlag (or caglad), their
paths are a.s. bounded on $[0,T]$, so
\begin{eqnarray} \label{boundednessCondsonCoeffs}
\limsup
_{h\rightarrow0}
\frac{\sup_{i} |\int_{t_{i-1}}^{t_i}a_s(\omega) \,\mathrm{d}s|}{h}
&\leq&
A(\omega)<\infty
 \quad \mbox{and}
 \nonumber
 \\[-8pt]
 \\[-8pt]
 \nonumber
\limsup _{h\rightarrow0} \frac{\sup_{i} |\int
_{t_{i-1}}^{t_i}\sigma^2_s(\omega)
\,\mathrm{d}s|}{h} &\leq&\Sigma(\omega)<\infty \qquad\mbox{a.s.}
\end{eqnarray}
It follows from
\cite{mancini07} that
\begin{equation} \label{DXzDivhlnunosuh}
\textup{a.s.}\qquad \sup_i \frac{|\int_{t_{i-1}}^{t_i}a_s\,\mathrm{d}s+\int
_{t_{i-1}}^{t_i}\sigma_s \,\mathrm{d}W_s|}{\sqrt{2h\log
{1}/{h} }}\leq
A\sqrt h+\sqrt\Sigma+1:=\Lambda.
\end{equation}
Since realistic values of $\sigma$
for asset prices belong to $[0.1,0.8]$ (in annual units), we
have that for small $h,$ the r.v. $\Lambda$ has order of
magnitude of 1, thus, in the finite jump intensity case, a.s. for
sufficiently small $h$, $(\Delta_{i} X)^2> r_h> 2h \log\frac{1}{h} $
indicates the presence of jumps in $]t_{i-1},t_{i}]$.

When $L$ has infinite activity, $\sum_{i=1}^n (\Delta_{i} X)^2I_{\{
(\Delta_{i} X)^2\leq r_h\}} $
behaves like $\sum_{i=1}^n (\Delta_{i} X)^2\times I_{\{ \Delta_{i} N=0,
|\Delta_{i}M|\leq
2\sqrt r_h\}}$ for small $h$ (Lemma \ref{IndDxqleqrvsIndDxdqleqr}).
Moreover, for any $\delta>0$, the jumps contributing to the
increments $\Delta_{i} X$ such that $(\Delta_{i} X)^2\leq r_h$ for
small $h$ have size
smaller than $c\sqrt{r_h + \delta}$ (\cite{mancini07}, Lemma 1), so
their contribution vanishes when $h\rightarrow0$. Note that $r_h=
ch^\beta$ satisfies condition \eqref{rh.eq} for any $\beta\in\,]0,1[$
and any constant $c$. Since
$\sqrt{2}\sigma\simeq1$ in most applications, we use $c=1$.
Define
\begin{equation} \eta^2(\varepsilon) := \int_{|x|\leq\varepsilon}
x^2 \nu(\mathrm{d}x),\qquad
d(\varepsilon) := \int_{\varepsilon< |x|\leq1} x \nu(\mathrm{d}x). \label
{cepsilon.eq}
\end{equation}
Let us remark that if $\lim _{h\rightarrow0}r_h= 0,$ then, by
\ref{asA1},
we have, as $h\rightarrow0$,
\begin{eqnarray} \label{OrdGrandIntPotxdnu}
\eta^2\bigl(2\sqrt{r_h}\bigr) &=& \int_{|x|\leq
2\sqrt{r_h}}
x^2 \nu(\mathrm{d}x)
\sim r_h^{1-{\alpha}/{2}},\qquad
\int_{|x|\leq2\sqrt{r_h}} x^k \nu(\mathrm{d}x) \sim
r_h^{{(k-\alpha)}/{2}},\nonumber\\
&&{}\hspace*{-40pt} k=3,4,\\
\nonumber
\int_{2\sqrt{r_h}<|x|\leq1} x \nu(\mathrm{d}x)&\sim&\bigl[c + r_h^{{(1-\alpha)}/{2}}\bigr]I_{\{\alpha\neq1\}} + \biggl[\ln\frac{1}{2\sqrt
{r_h}}
\biggr]I_{\{\alpha= 1\}} , \\
\int_{2\sqrt{r_h}<|x|\leq1} \nu(\mathrm{d}x)&\sim&
r_h^{-\alpha/2},\nonumber
\end{eqnarray}
where $\alpha$ is the BG index of $L$.
The following lemma, proved in the \hyperref[appendix.sec]{Appendix}, states that under~(\ref
{rh.eq}), each increment $\Delta_{i}M$ such that $|\Delta_{i}M|\leq
2\sqrt{r_h}$ only
contains jumps of magnitude less
than $2\sqrt{r_h}$ if $\alpha\leq1$, or smaller than $2h^{{1}/{(2\alpha)}} \log^{{1}/{(2\alpha)}}\frac1 h$ if $\alpha> 1$.

\begin{Lemma} \label{StructDXdwhenleqr} Define, for $h>0,$
$v_h:=h^{{1}/{(2\alpha)}} \log^{{1}/{(2\alpha)}}\frac1 h$. Under (\ref{rh.eq}).
there exists a sequence $h_k=T/n_k$ tending to zero as $k\rightarrow
\infty$
such that, for $k_0$ sufficiently large and $h\in\{ h_k,k\geq k_0 \}
$:
\begin{longlist}[(ii)]
\item[(i)] if $\alpha\leq1$, then for all $i=1,\ldots,n$,
\begin{eqnarray*}
&&\Delta_{i}M I_{\{ (\Delta_{i}M)^2\leq4 r_h\}}\\
&&\quad= 
\biggl(\int_{t_{i-1}}^{t_i}\!  \int_{|x|\leq2\sqrt{r_h}}   x
\tilde\mu(\mathrm{d}x, \mathrm{d}t) -
\int_{t_{i-1}}^{t_i}\!  \int_{ 2\sqrt{r_h}<|x|\leq1}   x
\nu(\mathrm{d}x)\,
\mathrm{d}t\biggr)
I_{\{ (\Delta_{i}M)^2\leq4 r_h\}}\qquad  a.s.;
\end{eqnarray*}
\item[(ii)] if $\alpha> 1$, then for all $i=1,\ldots,n$, we have
\begin{eqnarray*}
&&\Delta_{i}M I_{\{ (\Delta_{i}M)^2\leq4 r_h\}}\\
&&\quad= 
\biggl(\int_{t_{i-1}}^{t_i}\!  \int_{|x|\leq2 v_h}   x \tilde\mu
(\mathrm{d}x, \mathrm{d}t) -
\int_{t_{i-1}}^{t_i}\!  \int_{ 2 v_h<|x|\leq1}   x \nu(\mathrm{d}x) \,\mathrm{d}t
\biggr) I_{\{ (\Delta_{i}M)^2\leq4 r_h\}}
\qquad  a.s.
\end{eqnarray*}
%
\end{longlist}
\end{Lemma}

\begin{Remark}\label{JumpsLeqrOnefourth} Note that $v_h\leq r_h^{1/4}$
so that in the case \textup{(ii)} above ($\alpha>1$), for all $i=1,\ldots, n$, the
jumps of $M$ on $\{(\Delta_{i}M)^2\leq4r_h\}$ are bounded by $r_h^{1/4}$.
\end{Remark}
\begin{Definition*} Define
\begin{eqnarray} \label{defJh}
L^{(h)}_t&:=&\int_0^t\!\!\!\int_{|x|\leq 2\sqrt[4]r_h} x  \tilde\mu(\mathrm{d}x,
\mathrm{d}t) -
\int_0^t\!\!\! \int_{2\sqrt[4]r_h<|x|\leq1} x \nu(\mathrm{d}x) \,\mathrm{d}t,
\nonumber
\\[-8pt]
\\[-8pt]
\nonumber
 \Delta_{i} M^{(h)}&:=& \int_{t_{i-1}}^{t_i}\!\int_{|x|\leq2 \sqrt
[4]r_h} x
\tilde\mu(\mathrm{d}x, \mathrm{d}t).
\end{eqnarray}
\end{Definition*}
%
By Lemma \ref{StructDXdwhenleqr}, on a subsequence, a.s. for
sufficiently small $h$, $\forall i=1,\ldots, n$, on $\{(\Delta
_{i}M)^2\leq4r_h\}
,$ we
have
\begin{equation} \label{DefDecompDXd} \Delta_{i}M=\Delta_i L^{(h)}=
\Delta_{i} M^{(h)}- h d\bigl(2\sqrt[4]{ r_h}\bigr).
\end{equation}
$\Delta_{i} M^{(h)}$ is the compensated sum of jumps smaller in absolute
value than $2\sqrt[4]r_h$, while $h d(2\sqrt[4]{ r_h})$ is the
compensator of the
(missing) jumps larger than $2\sqrt[4]r_h$.

In \cite{mancini07}, a central limit theorem for $\hat{\mathit{IV}}_h$ was shown
in the case of
finite intensity jumps and cadlag adapted $\sigma$. Theorem
\ref{corVelCvIV} extends this to the case of infinite activity
without extra assumptions on $\sigma$. In particular, when
$\alpha<1$, the error $\hat{\mathit{IV}}_h -\mathit{IV}$ has the same rate of
convergence and asymptotic variance as in the case of finite
intensity jumps. The following proposition gives the asymptotic
variance of $(\hat{\mathit{IV}}_h-\mathit{IV})/\sqrt{2h}$ when $\alpha<1$.
\begin{Proposition}\label{teoStimIntQuart}
If $r_h=h^{\beta}$ with
$1> \beta> \frac{1}{2-\alpha/2}\in[1/2, 1[$, then, as $h\rightarrow0$,
\[
\hat{\mathit{IQ}}_h:=\frac{\sum_i (\Delta_{i} X)^4I_{\{ (\Delta_{i}
X)^2\leq r_h\}} }{3h} \stackrel{P}\rightarrow \mathit{IQ}= \int
_0^T \sigma_t^4 \,\mathrm{d}t.
\]
\end{Proposition}
%

The following result will be used to prove Theorem
\ref{corVelCvIV}. 
%
\begin{Theorem}\label{teoSpeedCvSumDXdq} Under Assumption \textup{\ref{asA1}},
as $h\to0,$
\begin{eqnarray}\label{SpeedCvSumDXdq}
&&\frac{\sum _{i=1}^n
(\int_{t_{i-1}}^{t_i}\!\int_{|x|\leq\varepsilon} x \tilde\mu
(\mathrm{d}x, \mathrm{d}t)- \int_{t_{i-1}}^{t_i}\!\int_{|x|\in
]\varepsilon, 1]} x \nu(\mathrm{d}x) \,\mathrm{d}t)^2
-T \ell_{2,h}\varepsilon^{2-\alpha}
-T\ell^2_{1,h} h \varepsilon^{2-2\alpha} I_{\{\alpha\neq1\}
}}{\sqrt T \sqrt
{\ell_{4,h}}\varepsilon^{2-\alpha/2}}
\nonumber
\\[-8pt]
\\[-8pt]
\nonumber
&&\quad\stackrel
{d}{\rightarrow} N(0,1),
\end{eqnarray}
where $\varepsilon=h^u$, $0<u\leq1/2$, $\ell_{j,h}= \int_{|x|\leq
\varepsilon} x^j \nu
(\mathrm{d}x) /  \varepsilon^{j-\alpha}$ for $j=2,4$ and
$\ell_{1,h}= \int_{\varepsilon<|x|\leq1} x \nu(\mathrm{d}x) /\break
[(c+\varepsilon^{1-\alpha
})I_{\{\alpha\neq1\}}+ \ln\frac{1}{2 \varepsilon} I_{\{\alpha=1\}
}]$ tend to
non-zero constants depending on $\nu$.
\end{Theorem}

We are now ready to state our central limit theorem for the
estimator $\hat{\mathit{IV}}_h$. A sequence $(X_{n})$ is said to converge
stably in law to a random variable $X$ (defined on an extension
$(\Omega',\mathcal{F}',P')$ of the original probability space) if
$\lim{E[Uf(X_{n})]}= E'[Uf(X)]$ for every bounded continuous
function $ f\dvtx \mathbb{R} \rightarrow\mathbb{R}$ and all bounded
random variables $ U$. This is obviously stronger than convergence
in law \cite{JacPro98}.
\vadjust{\eject}
\begin{Theorem}\label{corVelCvIV} Assume \textup{A1}
and $\sigma\not\equiv0$; choose $r_h=h^{\beta}$ with
$\beta> \frac{1}{2-\alpha/2} \in[1/2, 1[$. Then:
\begin{longlist}[(a)]
\item[(a)] if $\alpha<1$, 
we have, with $\stackrel{\mathit{st}}{\rightarrow}$ denoting stable
convergence in law,
\begin{equation} \label{hatIVasNor} \frac{\hat{\mathit{IV}_h}- \mathit{IV}}{\sqrt{2h
\hat{\mathit{IQ}}_h}}
\stackrel{\mathit{st}}\rightarrow{ N}(0, 1);
\end{equation}
\item[(b)]if $\alpha\geq1$, then
\[
\frac{\hat
{\mathit{IV}_h}- \mathit{IV}}{\sqrt{2h \hat{\mathit{IQ}}_h}} \stackrel{\mathit{a.s.}}\rightarrow+\infty.
\]
\end{longlist}
\end{Theorem}
\begin{Remark*} For $\alpha<1,$ Jacod \cite{jacod06}, Theorem 2.10(i),
has shown a
related central limit result for the threshold estimator of $\mathit{IV}$,
where $L$ is a semimartingale, but under the additional assumption
that $\sigma$ is an It\^o semimartingale. The proof of Theorem
\ref{corVelCvIV} in the case $\alpha<1$ does not rely on~\cite
{jacod06}, Theorem
2.10(i). An alternative proof under the It\^o
semimartingale assumption for $\sigma$ could combine the results
\cite{mancini07} with \cite{jacod06}, Theorem 2.10(i), in that
\begin{eqnarray*}
\frac{\hat{\mathit{IV}}-\mathit{IV}}{\sqrt h} &=& \frac{\hat{\mathit{IV}}(X_1)-\mathit{IV}}{\sqrt h} +
\frac
{\hat{\mathit{IV}}(M)}{\sqrt h}+
\frac{\sum_{i=1}^n (\Delta_{i} X_1)^2( I_{\{(\Delta_{i}
X)^2\leq r_h\}} - I_{\{(\Delta_{i} X_1)^2\leq
r_h\}})}{\sqrt{h}}\\
&&{}+
\frac{\sum_{i=1}^n (\Delta_{i}M)^2( I_{\{(\Delta_{i} X)^2\leq
r_h\}} - I_{\{(\Delta_{i}M)^2
\leq r_h\}})}{\sqrt{h}}\\
&&{}+ 2\frac{\sum_{i=1}^n \Delta_{i}
X_1\Delta_{i}M
I_{\{ (\Delta_{i} X)^2\leq r_h\}} }{\sqrt{h}},
\end{eqnarray*}
where
\[
\hat{\mathit{IV}}(X_1)\doteq\sum_{i=1}^n (\Delta_{i}
X_1)^2I_{\{(\Delta_{i} X_1)^2\leq
r_h\}},\qquad
\hat{\mathit{IV}}(M)\doteq
\sum_{i=1}^n (\Delta_{i}M)^2I_{\{(\Delta_{i}M)^2\leq r_h\}}.
\]
The first term
converges stably in law by \cite{mancini07}, the second one
converges stably to zero by \cite{jacod06}, Theorem 2.10(i). That
the remaining terms are negligible requires some further work (see the
proof of
Theorem \ref{corVelCvIV}).
\end{Remark*}

\section{Statistical tests}\label{testDescription.sec}

\subsection{Test for the presence of a continuous martingale component}
We now use the above results to design a test to detect the presence
of a continuous martingale component $\int_0^t\sigma_t \,\mathrm{d}W_t$, given
discretely recorded observations. Our test is feasible in the case
where $L$ has BG index $\alpha<1$, that is, the jumps are of finite
variation (see Section \ref{testAlpha.sec}). The test proceeds as
follows. First, we choose a coefficient $\beta\in[1/2, 1[$ close to
1. If we have an estimate $\hat{\alpha}$ of the BG index
\cite{WoernerHatAlpha06b,ASJ07,TodTau08}, then we may choose $\beta>
\frac{1}{2-\hat{\alpha}}$ (recall that $\frac{1}{2-\alpha}\in[1/2,
1[$). We choose a threshold $r_h=h^\beta$ and use the estimator $
\hat{\mathit{IQ}}_h$ of the integrated quarticity defined in Proposition
\ref{teoStimIntQuart}. We have shown in Theorem \ref{corVelCvIV}
that, when $\sigma\not\equiv0$ in the case $\alpha<1$, the
estimator $\hat{\mathit{IV}}_h$ is asymptotically Gaussian as $h\rightarrow0$.
However, if $\sigma\equiv0$, then both the numerator and the
denominator of (\ref{hatIVasNor}) tend to zero. To handle this case,
we
add an i.i.d. noise term:
\[
\Delta_{i} X^{v}:= \Delta_{i} X+v\sqrt{h} Z_i,\qquad  Z_i\mathop{\sim}^{\mathrm{i.i.d.}} {N}(0,1).
\]
As $h\rightarrow0,$
\[
\sum_{i=1}^n (\Delta_{i}
X^{v})^2\stackrel{P}{\rightarrow} [X^{v}]_T=
\int_0^T \sigma^2_s \,\mathrm{d}s +
v^2 T + T\int_{\mathbb{R}-\{0\}} x^2 \mu(\mathrm{d}x, \mathrm{d}s)
\]
and $ I_{\{ (\Delta_{i} X^{v})^2\leq r_h\}} $
removes the jumps of $X^{v}$ so that under the assumptions of
Theorem~\ref{corVelCvIV}, as $h\rightarrow0,$
\[
\hat{\mathit{IV}}^{v}_h := \sum_{i=1}^n (\Delta_{i} X^{v})^2I_{\{ (\Delta
_{i} X^{v})^2\leq r_h\}} \stackrel{P}\rightarrow\int
_0^T\sigma^2_s \,\mathrm{d}s + v^2 T.
\]
%
Under the null hypothesis\vspace*{-2pt} $\sigma\equiv0$,
we have $\hat{\mathit{IV}}^{v}_h \stackrel{P}{\rightarrow} v^2T$, $\hat{\mathit{IQ}}^{v}_h:=
\sum_i (\Delta_{i} X^{v})^4I_{\{ (\Delta_{i} X^{v})^2\leq r_h\}}
/\break (3h)\stackrel{P}{\rightarrow} v^4 T$
and
\begin{equation} \label{Threshclt} U_h:=\frac{\hat{\mathit{IV}}^{v}_h- v^2T}{
\sqrt{2h
\hat{\mathit{IQ}}^{v}_h}}\stackrel{\mathit{st}}\rightarrow\mathcal{N}.
\end{equation}
Note that if, on the contrary, $\sigma\not\equiv0$,
then we have that the limit in probability of $\hat{\mathit{IV}}^{v}_h$ 
is strictly larger than $ v^2 T$ and, by Lemma
\ref{IndDxqleqrvsIndDxdqleqr}, passing to a subsequence, a.s.
\begin{eqnarray*}
\lim _{h\rightarrow0}h \hat{\mathit{IQ}}^{v}_h &=&
\frac1 3 \lim _{h\rightarrow0}\sum_i (\Delta_{i}
X^{v})^4I_{\{ (\Delta_{i} X^{v})^2\leq r_h\}} =
\frac1 3 \lim _{h\rightarrow0} \sum_i (\Delta_{i}
X^{v})^4I_{\{ \Delta_{i} N=0, (\Delta_{i}M)^2\leq
2r_h\}}
\\
&\leq&\frac1 3 \lim _{h\rightarrow0} \sum_i \bigl(\Delta_{i}
X_0+\Delta_{i}M+v\sqrt{h}
Z_i \bigr)^4 I_{\{ (\Delta_{i}M)^2\leq2r_h\}}
\\
&\leq&\frac c 3 \lim _{h\rightarrow0} \sum_i (\Delta_{i}
X_0)^4 +
\frac c 3\lim _{h\rightarrow0}\sum_i (\Delta
_{i}M)^4I_{\{ (\Delta_{i}M)^2\leq2r_h\}}+
\frac c 3 \lim _{h\rightarrow0} \sum_i \bigl(v\sqrt{h} Z_i\bigr)^4.
\end{eqnarray*}
Using the facts that $\lim _{h\rightarrow0}\sum_i (\Delta
_{i}M)^4I_{\{
(\Delta_{i}M)^2\leq2r_h\}}\leq\lim _{h\rightarrow0}
2r_h\sum_i (\Delta_{i}M)^2I_{\{
(\Delta_{i}M)^2\leq2r_h\}}= 0$, by (\ref{VQThresholdedJumps}),
$\sum_i (\Delta_{i} X_0)^4/h \stackrel{P}{\rightarrow} c\int
_0^T\sigma^4_s\,\mathrm{d}s$ and
$\sum_i (v\sqrt{h} Z_i)^4/h \stackrel{\mathit{a.s.}}{\rightarrow} c v^4$, we
have, as
$h\rightarrow0$, $h \hat{\mathit{IQ}}^{v}_h \stackrel{P}{\rightarrow} 0.$
Therefore, under
the alternative $(H_1) \sigma\not\equiv0$, $U_h\to+\infty$ and
\mbox{$P\{|U_h|>1.96\} \rightarrow1$},
so the test is consistent.

\textit{Local power of the test}. To investigate the local power of
the test $U_h$, we consider a sequence of alternatives $(H_{1}^h)
\sigma=\sigma^h$, where $\sigma^h\downarrow0$. We denote by
$\hat{\mathit{IQ}}_{\sigma^h}^v, U_{\sigma^h}$ the statistics analogous to
$\hat{\mathit{IQ}}_{h}^v, U_h$, but constructed from $X^h_t=x_0+\int_0^t
a_s\,\mathrm{d}s+\int_0^t\sigma^h_s \,\mathrm{d}W_s+ L_t,   t\in\,]0,T].$ In the case of
constant $\sigma$ and $\sigma^h$, and finite jump intensity, using
standard results on convergence of sums of a triangular array
\cite{JacDisp07}, Lemmas 4.1 and 4.3, we have
\[
\hat{\mathit{IQ}}_{\sigma^h}^v\stackrel{\mathrm{ucp}}{\rightarrow} v^4T,\qquad
U_{\sigma^h}
\stackrel{d}{\rightarrow}
\lim_{h\rightarrow0} \frac{(\sigma^h)^2 }{\sqrt h} T + \sqrt2 v^2 Z_T,
\]
%
where $\stackrel{\mathrm{ucp}}{\rightarrow}$ denotes uniform convergence in
probability on compacts subsets of $[0,T]$ \cite{Protter05} and $Z$
is a standard Brownian motion. So, either $U_{\sigma^h}$ tends in
distribution to $c+\sqrt2 v^2 Z_T$, if $\sigma^h =\mathrm{O}(h^{1/4}),$ or
$U_{\sigma^h}\rightarrow\infty$, if $h^{1/4}=\mathrm{o}(\sigma^h)$. Thus, if
$c$ is a
(possibly zero) constant, we have:
\begin{eqnarray*}
&&\mbox{if } \frac{\sigma^h}{h^{1/4}} \rightarrow c,\qquad \mbox{then } P\{
U_{\sigma
^h}>1.64|H_1^h\}
\rightarrow P\biggl\{Z_1>\frac{1.64-c^2T}{\sqrt{2T} v^2}\biggr\};
\\
&&\mbox{if } \frac{\sigma^h}{h^{1/4}} \rightarrow+\infty,\qquad \mbox{then } P\{
U_{\sigma^h}>1.64|H_1^h\}
\rightarrow1.
\end{eqnarray*}
For values of $v$ in Section \ref{simuls.sec}, we have
$1.64/\sqrt{2T} v^2=\mathrm{O}(10^8)$ and thus the local power of the test is
small if $\sigma^h =\mathrm{O}(h^{1/4})$.

\subsection{Testing whether the jump component has finite
variation}\label{testAlpha.sec}

To construct a test for discriminating $\alpha<1$ from $\alpha
\geq1$, Theorem \ref{corVelCvIV} suggests the use~of $(\hat{\mathit{IV}}_h -
\mathit{IV})/\sqrt{2h\hat{\mathit{IQ}}_h}$, but this requires knowing the process
$\sigma$ to compute $\mathit{IV}$. We propo\-se a feasible alternative.
Consider, instead, the estimator
\[
\hat H_{h}:= \sum_{i=1}^n \Delta_i
X I_{\{(\Delta_{i} X)^2> r_h\}}= X_T-\sum_{i=1}^n \Delta_i X I_{\{
(\Delta_{i} X)^2\leq
r_h\}}.
\]
\begin{Proposition}\label{EstimJwhenAlphaLess1}
When $\alpha<1$, $\hat H_{h }$ is a consistent estimator of $J_T + m T
$, $m:= \int_{-1}^1 x \nu(\mathrm{d}x)$.
\end{Proposition}

Consider $Z_i=\Delta_i W^v$, where $W^v$ is a Wiener process
independent of $W,L,$
and define
\[
\Delta_i\hat H^{v}:= \Delta_i X I_{\{(\Delta_{i} X)^2> r_h\}}
+v\sqrt{h} Z_i
\quad \mbox{and}\quad   H^v_T:= J_T + mT + vW^v_T.
\]
Under the null hypothesis $\alpha<1,$
\[
\hat{\mathit{IV}}_h^{H^v}:= \sum_i (\Delta_i\hat H^{v})^2 I_{\{(\Delta_i\hat
H^{v})^2 \leq r_h\}}
\]
is an estimator of the integrated variance $v^2T$ of
$H^{v}$, so, under the null hypothesis $\mathrm{(H_0)}$ $\alpha<1,$
we can find $\beta>\frac{1}{2-\alpha}\in\,]\frac1 2,1[$ such that
\begin{equation} \label{Uhalpha}U^{(\alpha)}_h:= \frac{\hat
{\mathit{IV}}_h^{H^v} - v^2T}{
\sqrt{2h\hat{\mathit{IQ}}_h^{H^v}}} \stackrel{d}\rightarrow{N}(0,1),
\end{equation}
where $\hat{\mathit{IQ}}_h^{H^v}:= \frac{1}{3h}\sum_i (\Delta_i \hat H
^{v})^4 I_{\{(\Delta_i \hat H^{v})^2\leq r_h\}}$ and
$r_h=h^\beta$. In particular, \mbox{$ P\{|U^{(\alpha)}_h|>1.96\}$}$\rightarrow5\%$.

If, on the contrary, $\alpha\geq1$, then reasoning as in Theorem \ref
{corVelCvIV},
for any $\beta\in\,]0,1[$, we have $
U^{(\alpha)}_h \stackrel{P}\rightarrow+\infty,$ so the test is
consistent. If
$|U^{(\alpha)}_h|>1.96$, then we reject $\mathrm{(H_0)}$ $\alpha<1$ at the 95\% confidence
level.

\begin{Remark*} To apply this test, we first need to
decide whether $\alpha< 1$, using the previously described test.
\end{Remark*}

\begin{table}[b]
\tablewidth=285pt
\caption{Testing for finite variation of jumps: $\alpha$-stable process plus
Brownian motion. $\mathit{pct}$ is the percentage of outcomes where
$|U^{(\alpha)}_{(j)}|> 1.96$}\label{tab1}
\begin{tabular*}{285pt}{@{\extracolsep{\fill}}lld{1.6}@{\qquad}ll@{\qquad}ll@{}}
\hline
$n$& $h$ &\multicolumn{1}{l}{$v$} &$\alpha$& $\mathit{pct}$&$\alpha$ & $\mathit{pct}$\\
\hline
\phantom{0\,}1000 & 5 min & 0.000001 &0.6 & 0.067 &1.6 & 0.439\\
\phantom{0\,}1000 & 5 min & 0.0001 &0.6 & 0.056 &1.6 & 0.407\\
\phantom{0\,}1000 & 5 min & 0.01 &0.6 & 0.047 &1.6 & 0.250\\
\phantom{0\,}1000 & 5 min & 0.1 &0.6 & 0.053 &1.6 & 0.726\\
\phantom{0\,}1000 & 1 min & 0.0001 &0.6 & 0.049 &1.6 & 0.241\\
\phantom{0\,}1000 & 1 hour& 0.0001 &0.6 & 0.051 &1.6 & 0.875\\
\phantom{0\,}1000 & 1 day & 0.0001 &0.6 & 0.066 &1.6 & 0.984\\
\phantom{00\,}100 & 5 min & 0.0001 &0.6 & 0.065 &1.6 & 0.137\\
10\,000 & 5 min & 0.0001 &0.6 & 0.065 &1.6 &
0.928\\
\hline
\end{tabular*}
\end{table}

\section{Numerical experiments}\label{simuls.sec}

\subsection{Testing the finite variation of the jump component}\label
{sectTestforFV}

We simulate $n$ increments $\Delta_{i} X$ of a process $ X= \sigma W + L,$
where $L$ is a symmetric $\alpha$-stable L\'evy process, 
$\sigma=0.2$. We generate 1000 independent samples containing $n$
increments each and compute $U^{(\alpha)}_h$ as in (\ref{Uhalpha})
for a range of values of $v$, $h$ (1 minute, 5 minutes, 1 hour, 1
day) and number of observations $n$. Table \ref{tab1} reports the
percentage (\textit{pct}) of outcomes where
$|U^{(\alpha)}_{h(j)}|
>1.96$, $j=1,\ldots, 1000$,
for threshold exponent $\beta
=0.999$. Note that with $n=1000$ and $h$ equal to five minutes
($h=1/(252\times84)$), we have $T <1$ year; for $\alpha=0.6,$ the
lower bound for $\beta$ is $\frac{1}{2-\alpha}=0.71$; when $n=1000,
h=1/(84\times252)$ and the BG index of $L$ is 0.6 (resp., 1.6), the
ratio of $v=10^{-4}$ to the standard deviation of the increments
$\Delta_{i} X$
is 0.074 (resp., 0.022).
The test results
are observed to be reliable if we use $n=10\thinspace000$ observations,
a time
resolution of five minutes and $v=10^{-4}$. In fact, when the data-generating
process has BG index 0.6, the test leads us to accept
the hypothesis $\mathrm{(H_0)}$ $\alpha<1$ in about 94 cases out of 100. On the
contrary, when the
process has BG index 1.6, the test tells us to reject $\mathrm{(H_0)}$ in 92
cases out of 100.


%
%


\subsection{Test for the presence of a Brownian component}\label
{SectionTestPresSig}
We simulate 1000 independent paths 
of a process $ X_t = \int_0^t \sigma_u \,\mathrm{d}W_u + L,$
for different L\'evy processes $L$ and constant or stochastic
$\sigma$, on a time grid with $n$ steps. We take threshold
$r_h=h^{0.999}$. For each trial $j=1,\ldots,1000$,
we compute $U_{h(j)}$ given in (\ref{Threshclt}) and report the
percentage (\textit{pct}) of cases where $|U_{h(j)}|>1.96$.

\begin{example}[(Brownian motion
plus compound Poisson
process, BG index $\bolds{\alpha=0}$)]\label{BMplusCPTestPresSigma} We consider here constant
$\sigma$ and $L=\sum_{i=1}^{N_t} B_i$, a compound Poisson process
with i.i.d. ${N}(0, 0.6^2)$ sizes of jump and jump intensity
$\lambda=5$ (as in \cite{aitsahalia04}). Table \ref{tab2} illustrates
the performance of our test for various time steps $h$, numbers of
observations $n$ and noise levels $v$:
%
\begin{table}[b]
\tablewidth=285pt
\caption{Testing for the presence of a Brownian component: case of
Brownian motion plus compound Poisson jumps (Example \protect\ref{BMplusCPTestPresSigma})}\label{tab2}
\begin{tabular*}{285pt}{@{\extracolsep{\fill}}lld{1.6}@{\qquad}ll@{\qquad}ll@{}}
\hline
$n$& $h$ &\multicolumn{1}{l}{$v$} &$\sigma$& $\mathit{pct}$&$\sigma$ & $\mathit{pct}$\\
\hline
\phantom{0\,}1000 & 5 min & 0.000001 &0 & 0.043 &0.2 & 1\\
\phantom{0\,}1000 & 5 min & 0.0001 &0 & 0.048 &0.2 & 1\\
\phantom{0\,}1000 & 5 min & 0.01 &0 & 0.054 &0.2 & 1\\
\phantom{0\,}1000 & 5 min & 0.1 &0 & 0.041 &0.2 & 1\\
\phantom{0\,}1000 & 1 min & 0.0001 &0 & 0.047 &0.2 & 1\\
\phantom{0\,}1000 & 1 hour& 0.0001 &0 & 0.054 &0.2 & 1\\
\phantom{0\,}1000 & 1 day & 0.0001 &0 & 0.082 &0.2 & 1\\
\phantom{00\,}100 & 5 min & 0.0001 &0 & 0.065 &0.2 & 1\\
10\,000 & 5 min & 0.0001 &0 & 0.049 &0.2 &
1\\
\hline
\end{tabular*}
\end{table}
Note that when $\sigma=0$ (resp.,~0.2), $n=1000$ and $h=1/(84\times
252)$ the ratio of $v=10^{-4}$ to the standard deviation of the returns
$\Delta_{i} X$ equals 0.007 (resp., 0.052).
\end{example}

We find that the test is reliable for
values $n=1000$, $h=5$ minutes and $v=10^{-4}$ since it correctly
accepts $\mathrm{(H_0)}$ in 95 cases out of 100 and rejects $\mathrm{(H_0)}$ in all
cases when it is false.

\begin{example}[({Brownian motion plus $\alpha$-stable jumps: $\alpha\in\,]0,2[$})]\label{ex2}
Here, $L$ is a~symmetric $\alpha$-stable L\'evy process and
$\sigma$ is constant. The results in Table \ref{tab3} confirm the satisfactory
performance of the test when $\alpha=0.3 <1$ for $n=1000$, $h=5$
minutes and $v=10^{-4}$.
%
\begin{table}
\tablewidth=285pt
\caption{Testing for the presence of a Brownian component: case of
Brownian motion plus $\alpha$-stable L\'evy process with $\alpha=0.3$
(Example \protect\ref{ex2})}\label{tab3}
\begin{tabular*}{285pt}{@{\extracolsep{\fill}}lld{1.6}@{\qquad}ll@{\qquad}ll@{}}
\hline
$n$& $h$ &\multicolumn{1}{l}{$v$} &$\sigma$& $\mathit{pct}$&$\sigma$ & $\mathit{pct}$\\
\hline
\phantom{0\,}1000 & 5 min & 0.000001 &0 & 0.042 &0.2 & 1\\
\phantom{0\,}1000 & 5 min & 0.0001 &0 & 0.026 &0.2 & 1\\
\phantom{0\,}1000 & 5 min & 0.01 &0 & 0.054 &0.2 & 1\\
\phantom{0\,}1000 & 5 min & 0.1 &0 & 0.053 &0.2 & 1\\
\phantom{0\,}1000 & 1 min & 0.0001 &0 & 0.046 &0.2 & 1\\
\phantom{0\,}1000 & 1 hour& 0.0001 &0 & 0.140 &0.2 & 1\\
\phantom{0\,}1000 & 1 day & 0.0001 &0 & 0.805 &0.2 & 1\\
\phantom{00\,}100 & 5 min & 0.0001 &0 & 0.056 &0.2 & 1\\
10\,000 & 5 min & 0.0001 &0 & 0.165 &0.2 &
1\\
\hline
\end{tabular*}
\end{table}
%

Table \ref{tab4}, for the case $\alpha=1.2>1$, confirms that
we cannot rely on the test results in this case: even when $\sigma
\equiv0$, the statistic $U_h$ diverges if $\alpha\geq1$.

\begin{table}[b]
\tablewidth=285pt
\caption{Testing for the presence of a Brownian component: case of
Brownian motion plus $\alpha$-stable L\'evy process with $\alpha=1.2$
(Example \protect\ref{ex2})}\label{tab4}
\begin{tabular*}{285pt}{@{\extracolsep{\fill}}lld{1.6}@{\qquad}ll@{\qquad}ll@{}}
\hline
$n$& $h$ &\multicolumn{1}{l}{$v$} &$\sigma$& $\mathit{pct}$&$\sigma$ & $\mathit{pct}$\\
\hline
\phantom{0\,}1000 & 5 min & 0.000001 &0 &1 &0.2 & 1\\
\phantom{0\,}1000 & 5 min & 0.0001 &0 &1 &0.2 & 1\\
\phantom{0\,}1000 & 5 min & 0.01 &0 &1 &0.2 & 1\\
\phantom{0\,}1000 & 5 min & 0.1 &0 &1 &0.2 & 1\\
\phantom{0\,}1000 & 1 min & 0.0001 &0 &1 &0.2 & 1\\
\phantom{0\,}1000 & 1 hour& 0.0001 &0 &1 &0.2 & 1\\
\phantom{0\,}1000 & 1 day & 0.0001 &0 &1 &0.2 & 1\\
\phantom{00\,}100 & 5 min & 0.0001 &0 &0.994 &0.2 & 1\\
10\,000 & 5 min & 0.0001 &0 &1 &0.2 & 1\\
\hline
\end{tabular*}
\end{table}\vspace*{-3pt}
\end{example}

The main point here is that we may use a model-free
choice of threshold.

\begin{example}[(Stochastic volatility plus variance gamma
jumps: $\bolds{\alpha=0}$)]\label{ex3} Let us now consider a model $X$ with
stochastic volatility $\sigma_t$, correlated with the Brownian motion
driving $X$ and with jumps given by an independent variance gamma process:
%
\[
\mathrm{d}X_t= (\mu-\sigma_t^2/2) \,\mathrm{d}t+\sigma_t \,\mathrm{d}W_t^{(1)}+ \mathrm{d}L_t,
\]
where
\begin{equation} \label{dynamicSigma} \sigma_t = \mathrm{e}^{K_t},\qquad  \mathrm{d} K_t= -k
(K_t - \bar K) \,\mathrm{d}t + \varsigma \,\mathrm{d}W_t^{(2)},\qquad  d\bigl\langle W^{(1)}, W^{(2)}\bigr\rangle_t
= \rho \,\mathrm{d}t,
\end{equation}
$W^{(\ell)}$ are standard Brownian motions, $\ell=1,2,3,$ and
$L_t=c G_t + \eta W^{(3)}_{G_t}$ is an independent variance gamma
process, a pure-jump L\'evy process with
BG index $\alpha=0$ \cite{Madan01}; $G$ is a gamma subordinator
independent of $W^{(3)}$ with $G_{h}\sim\Gamma(h/b, b)$. For
$\sigma,$ we choose $K_0= \ln(0.3)$, $k=0.09$, $\bar K=\ln(0.25)$,
$\varsigma=0.05$ to ensure that $\sigma$ fluctuates in the range
0.2--0.4. As for the jump part of $X,$ we use $\operatorname{Var}(G_1)=b= 0.23$,
$\eta=0.2$, $c=-0.2$, estimated from the S\&P 500 index in~\cite{Madan01}. The remaining parameters are $\rho=-0.7$ and
$\mu=0$. The following results in Table~\ref{tab5} confirm the reliability of the test
for the presence of a Brownian component with $n=1000$, $h=5$ minutes
and $v=10^{-4}$.
%
\begin{table}[b]\vspace*{-3pt}
\tablewidth=285pt
\caption{Testing for the presence of a Brownian component: stochastic
volatility process with variance gamma jumps (Example \protect\ref{ex3})}\label{tab5}
\begin{tabular*}{285pt}{@{\extracolsep{\fill}}lld{1.6}@{\qquad}ll@{\qquad}ll@{}}
\hline
$n$& $h$ &\multicolumn{1}{l}{$v$} &$\sigma$& $\mathit{pct}$&$\sigma$ & $\mathit{pct}$\\
\hline
\phantom{0\,}1000 & 5 min & 0.000001 &0 & 0.032 &Stoch. & 1\\
\phantom{0\,}1000 & 5 min & 0.0001 &0 & 0.017 &Stoch. & 1\\
\phantom{0\,}1000 & 5 min & 0.01 &0 & 0.027 &Stoch. & 1\\
\phantom{0\,}1000 & 5 min & 0.1 &0 & 0.054 &Stoch. & 1\\
\phantom{0\,}1000 & 1 min & 0.0001 &0 & 0.034 &Stoch. & 1\\
\phantom{0\,}1000 & 1 hour& 0.0001 &0 & 0.918 &Stoch. & 1\\
\phantom{0\,}1000 & 1 day & 0.0001 &0 & 1.000 &Stoch. & 1\\
\phantom{00\,}100 & 5 min & 0.0001 &0 & 0.049&Stoch. & 1\\
10\,000 & 5 min & 0.0001 &0 & 0.912 &Stoch.
& 1\\
\hline
\end{tabular*}
\end{table}
\end{example}

\begin{Remark*} In \cite{ManRen08}, a variable threshold
function is used to estimate the volatility, in order to account for
heteroscedasticity and volatility clustering, with\vadjust{\goodbreak}
results
very similar to the ones obtained with a constant threshold. This is
justified by the fact
that in most applications, values of~$\sigma$ are within the range
[0.1, 0.8], thus the order of magnitude of $\Lambda$ in (\ref
{DXzDivhlnunosuh}) is o 1.
\end{Remark*}

\section{Applications to financial time series}\label{data.sec}

We apply our tests to explore the fine structure of price
fluctuations in two financial time series. We consider the DM/USD
exchange rate from October 1st, 1991 to November 29th, 1994 and
the SPX futures prices from January 3rd, 1994 to December 18th, 1997. From
high-frequency time series, we build five-minute log-returns
(excluding, in the case of SPX futures, overnight log-returns). This
sampling frequency avoids many microstructure effects seen at
shorter time scales (e.g., seconds), while leaving us with a
relatively large sample.

\subsection{Deutsche Mark/USD exchange rate}
\begin{figure}[b]

\includegraphics{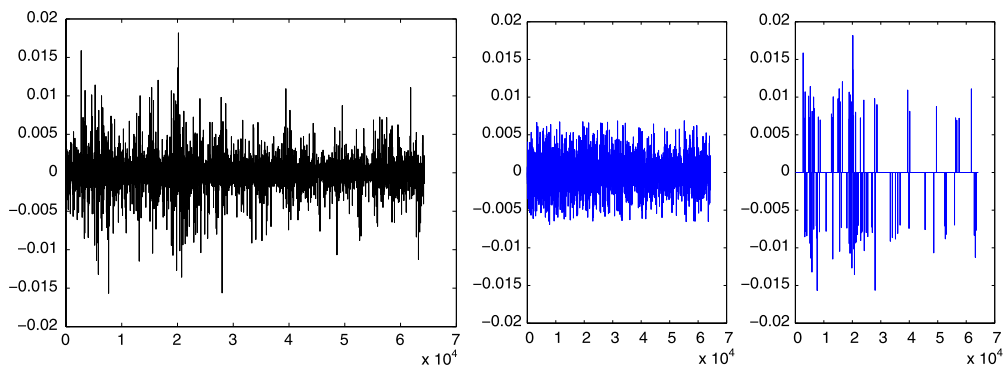}

  \caption{Left: DM/USD five-minute log-returns, October 1991 to November
1994. Center: plot of
$\Delta_{i} XI_{\{ (\Delta_{i} X)^2\leq r_h\}}$, $i=1,\ldots,n$.
Right: increments with jumps $\Delta_{i} XI_{\{ (\Delta_{i} X)^2>
r_h\}}
, i=1,\ldots,n$.}\label{plotDiffLogDM}
\end{figure}

The DM/USD exchange
rate time series was compiled by Olsen \& Associates. We consider the
series of 64\thinspace284 equally spaced five-minute log-returns, with
$h=\frac{1}{252\times84}\approx4.7 \times
10^{-5},$ displayed in Figure \ref{plotDiffLogDM}.
%

Barndorff-Nielsen and Shephard \cite{bnsBPV06a} provide evidence
for the presence of jumps in this series using nonparametric
methods. Using as threshold $r_h=h^{0.999}$, we apply the test of
Section \ref{sectTestforFV} to the degree of activity of the jump
component. As in the simulation study, we divide the data into 64
non-overlapping batches of $n=1000$ observations each and compute,
for each batch, the statistic $U_{h(j)}^{(\alpha)}$, $j=1,\ldots,64$, with
$v= 10^{-4}$. Only $4.7\%$ of the values observed are outside the
interval $[-1.96, 1.96]$, hence we cannot reject the assumption
$\mathrm{(H_0)}$ $\alpha<1$. Given this result, we can now use the test in
Section \ref{SectionTestPresSig} for
the presence of a Brownian component in the price process.
Computation of the statistic $U_h$ shows values much larger
than 1.96 for all batches: we reject $\mathrm{(H_0)}$ $\sigma\equiv0$.
These results indicate, for instance, that a variance gamma model,
with no Brownian component, would be inadequate for the DM/USD time
series.

\subsection{S\&P 500 index}
\begin{figure}

\includegraphics{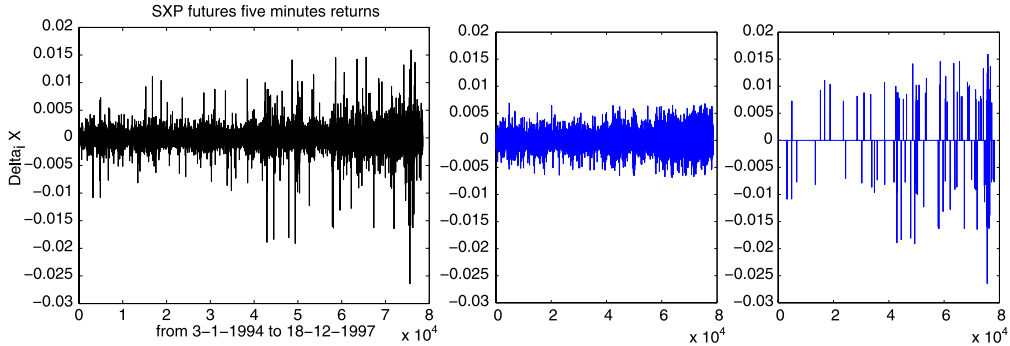}

  \caption{Left: SPX five-minute log-returns, January 1994 to December 1997.
Center: plot of
$\Delta_{i} XI_{\{ (\Delta_{i} X)^2\leq r_h\}}$, $i=1,\ldots, n$.
Right: increments with jumps $\Delta_{i} XI_{\{ (\Delta_{i} X)^2>
r_h\}}
, i=1,\ldots,n$.}\label{plotDiffLogSPX}
\end{figure}

We consider a series of 78\,497 non-overlapping five-minute
log-returns, as displayed in Figure \ref{plotDiffLogSPX}. Using as
threshold $r_h=h^{0.999}$, we decompose the series into periods
displaying jumps and other periods, as displayed in Figure \ref
{plotDiffLogSPX} (central and right panels).

We divide the data into 78
non-overlapping batches of $n=1000$ observations each and compute, for
each batch,
the statistic $U_{h(j)}^{(\alpha)}$, $j=1,\ldots,64$, with $v= 10^{-4}$.
$5.1\%$ of the values observed are
outside the interval $[-1.96, 1.96]$: for this period, we cannot reject
the assumption $\mathrm{(H_0)}$ $\alpha<1$.
Given this result, we can use the test for
the presence of a Brownian component in the price process.
Computation of the statistic $U_h$ shows values much larger
than 1.96 for all batches: we reject $\mathrm{(H_0)}$ $\sigma\equiv0$. The test
thus indicates
the presence of a Brownian martingale component.

We note that our findings contradict the conclusion of Carr \textit{et
al.} \cite{cgmy02} who model the (log-) SPX index from 1994 to 1998 as
a tempered stable L\'evy process plus a Brownian
motion and propose
a pure-jump model using a parametric estimation method. Under less
restrictive assumptions on the structure of the process and using
our nonparametric test, we find evidence for a non-zero Brownian
component in
the index.

\begin{appendix}
\section*{Appendix: Technical results and proofs}\label{appendix.sec}

\begin{pf*}{Proof of Lemma \protect\ref{StructDXdwhenleqr}}
By \cite{Met82}, Theorem 25.1, there exists a sequence $(n_k)$
such that
\renewcommand{\theequation}{\arabic{equation}}
\setcounter{equation}{16}
\begin{equation} \label{CvQvUcpMetiv} \sup_{t_j\in\Pi^{(n_k)}}
\biggl|(\Delta_{j}M)^2- \sum_{s\in]t_{j-1}, t_j]} (\Delta
M_{s})^2\biggr|\stackrel{\mathrm{a.s.}}{\rightarrow} 0,
\end{equation}
where $\Pi^{(n_k)}$ is the
partition of $[0,T]$ on which the increments $(\Delta_{i}M)^2$ are
constructed. Let us rename $n_k$ as $n$. Using It\^o's formula,
we have
\[
(\Delta_{i}M)^2- \sum_{s\in]t_{i-1}, t_i]} (\Delta M_{s})^2=
2\int_{t_{i-1}}^{t_i}( M_{ s_{-}}- M_{t_{i-1}} )\,\mathrm{d}M_s.
\]
\begin{longlist}[(ii)]
\item[(i)] For $\alpha< 1,$ our statement is proved in \cite{ManRen08}, Lemma
A.2, which uses the fact that the speed of convergence to 0 of $\sum
_{i=1}^n |\int_{t_{i-1}}^{t_i}(M_{ s_{-}}- M_{ t_{i-1}}) \,\mathrm{d}M_{ s}|$ is
shown in \cite
{jacod04} to be $u_n= n$. For $\alpha=1,$ the same reasoning
can be repeated since $u_n=n/(\log n)^2$ does not change the
conclusion.

\item[(ii)] If $\alpha>1,$ we have $u_n= (n/\log n)^{1/\alpha}$ and can only
conclude that a.s. for small $h,$
%
\[
\sup_i \biggl|\int_{t_{i-1}}^{t_i}(M_{ s_{-}}- M_{ t_{i-1}}) \,\mathrm{d}M_{
s}\biggr|\leq c u_n^{-1}
\]
with $c>0$, so that a.s. for small $h,$ we have
\begin{eqnarray*}
\sup_i \biggl(\sum_{s\in]t_{i-1},t_i]} (\Delta M_{s})^2\biggr) I_{\{
(\Delta_{i}M)^2\leq4r_h\}
}&\leq&\sup_i \biggl|(\Delta_{i}M)^2- \sum_{s\in]t_{i-1},t_i]}
(\Delta M_{s})^2\biggr|+\sup_i
|(\Delta_{i}M)^2|
\\
%
%
&\leq& c u_n^{-1}+4r_h=\mathrm{O}\biggl(\delta^{1/ \alpha} \log^{1/ \alpha
}\frac1 h\biggr).
\end{eqnarray*}
\end{longlist}
\upqed\end{pf*}
\begin{Lemma} \label{PDX0grandeePDtildeJ2grande0} 
Under \textup{(\ref{rh.eq})}:
\begin{longlist}[(ii)]
\item[(i)] there exists a strictly positive variable $\bar{h}$ such that
for all $ i=1,\ldots,n,$
\renewcommand{\theequation}{\arabic{equation}}
\setcounter{equation}{17}
\begin{equation} \label{DXzeromagr}
I_{\{h\leq\bar{h}\}} I_{\{(\Delta_{i}X_0)^{2}>r_h\}}=0 \qquad  \mbox{a.s.};%
\end{equation}
\item[(ii)]
\begin{equation}\label{DNneq0eDXdgeqr}
\forall c>0,  n
P\{\Delta_{i} N\neq0,(\Delta_{i}M)^2> cr_h\}\mathop{\longrightarrow}^{h\to0} 0;
\end{equation}
\item[(iii)] in
the case $r_h=h^\beta$, $\beta\in\,]0,1[$, we have
\begin{equation}
\label{SumIDXqmagr}\limsup_{h\rightarrow
0}h^{{\alpha\beta}/{2}}\sum_{i=1}^n P\{(\Delta_{i} X)^2> r_h\}
\leq c.
\end{equation}
\end{longlist}
\end{Lemma}
\begin{pf}
Equality (\ref{DXzeromagr}) is a consequence of
(\ref{DXzDivhlnunosuh}), while (\ref{DNneq0eDXdgeqr}) is a
consequence of the independence of $N$ and $M$, and of the
Chebyshev inequality: as $h\rightarrow0,$
\[
n P\{\Delta_{i} N\neq0, (\Delta_{i}M)^2> cr_h\} \leq n \mathrm{O}(h)\cdot
\frac{E[(\Delta_{i}M)^2]}{cr_h
}= \mathrm{O}\biggl(\frac{h}{r_h}\biggr).
\]
The proof of (\ref{SumIDXqmagr}) can be achieved as in \cite{ASJ07},
Lemma 6, but
we give a simpler proof under our assumptions.
It is sufficient to show that
\renewcommand{\theequation}{\arabic{equation}}
\setcounter{equation}{20}
\begin{equation} \label{SumIDXqmagrLevy}
P\{(\Delta_{i} X)^2> r_h\}\leq c h^{1-{\alpha\beta}/{2}}.
\end{equation}
First, 
we show that
\begin{equation} \label{DXmagSogliaRidottoADXd}
P\bigl\{|\Delta_{i} X| >\sqrt{r_h}\bigr\} = P\bigl\{|\Delta_{i}M| >\sqrt{r_h}/4\bigr\}+
\mathrm{O}(h^{1-\alpha\beta/2})
\end{equation}
so that for (\ref{SumIDXqmagrLevy}), it is sufficient to prove that
\begin{equation} \label{PDXdmagSoglia}
P\bigl\{|\Delta_{i}M| > \sqrt{r_h}/4\bigr\}\leq c h^{1-{\alpha\beta}/{2}}.
\end{equation}
%
To show (\ref{DXmagSogliaRidottoADXd}), note that
if $|\Delta_{i} X| >\sqrt{r_h}$, then either $\Delta_{i} J\neq0$ or
$|\Delta_{i}M| >\sqrt{r_h}/4$ since,
for small~$h$,
\begin{equation} \label{doveDXqgrr} \sqrt{r_h}< |\Delta_{i} X|\leq
|\Delta_{i} X_0|+ |\Delta_{i} J|+ |\Delta_{i}M| \leq\sqrt{r_h}/2 +
|\Delta_{i} J|+ |\Delta_{i}M|\qquad
\mbox{a.s.}
\end{equation}
Thus,
\[
P\bigl\{|\Delta_{i} X| >\sqrt{r_h}\bigr\} \leq P\{\Delta
_{i} J\neq0 \} + P\bigl\{|\Delta_{i}M|
>\sqrt{r_h}/4\bigr\}
\]
and since $P\{\Delta_{i} J\neq0 \}=\mathrm{O}(h)=\mathrm{o}(h^{1-\alpha\beta/2}),$
(\ref
{DXmagSogliaRidottoADXd}) is verified.

%

In order to verify (\ref{PDXdmagSoglia}), define
$\tilde N_t :=\sum_{s\leq t} I_{\{|\Delta M_{s} |>\sqrt{r_h}/4\}}$
and write
\begin{eqnarray}\label{ProbDtildeNgeq1}
P\bigl\{|\Delta_{i}M|>\sqrt{r_h}/4\bigr\} &=& P\bigl\{\Delta_{i} \tilde N=0, |\Delta
_{i}M|> \sqrt{r_h}/4\bigr\}\nonumber\\ 
&&{}+ P\bigl\{\Delta_{i} \tilde N\geq1, |\Delta_{i}M|>
\sqrt{r_h}/4\bigr\}\\
&\leq& P\{\Delta_{i} \tilde
N\geq1\} + 
P\bigl\{\Delta_{i} \tilde N=0, |\Delta_{i}M|>
\sqrt{r_h}/4\bigr\}.\nonumber
\end{eqnarray}
Note that $ \tilde N_t
=\int_0^t\!\!\!\int_{|x|>\sqrt{r_h}/4} \mu(\mathrm{d}x, \mathrm{d}t)$ is a compound Poisson
process with intensity
$\nu\{|x|>\sqrt{r_h}/4\}= \mathrm{O}(r_h^{-\alpha/2})$, so
$ P\{\Delta_{i} \tilde N\geq1\}=\mathrm{O}(h\nu\{|x|>\sqrt{r_h}/4\})=
\mathrm{O}(h^{1-\alpha\beta/2})$
and thus the first term above is dominated by $h^{1-\alpha\beta/2}$,
as required.
Finally, on $\{\Delta_{i} \tilde N=0\}$, $M$ does not have jumps
bigger than $\sqrt{r_h}/4$
on the interval $]t_{i-1}, t_i]$, so
\[
\Delta_{i}M= \int_{t_{i-1}}^{t_i}\int_{|x|\leq\sqrt{r_h}/4} x
\tilde\mu(\mathrm{d}x, \mathrm{d}t)- h \int_{\sqrt{r_h}
/4<|x|\leq1} x \nu(\mathrm{d}x),
\]
therefore
\begin{eqnarray*}
P\bigl\{\Delta_{i} \tilde N=0, |\Delta_{i}M|> \sqrt{r_h}/4\bigr\}
&\leq& P\bigl\{|\Delta_{i}M|> \sqrt{r_h}/4, |\Delta M_{s} |\leq\sqrt
{r_h}/4 \mbox{ for all
} s\in\,]t_{i-1},t_i]\bigr\}
\\
&\leq&4\frac{E[(\Delta_{i}M)^2I_{\{ |\Delta M_{s} |\leq\sqrt
{r_h}/4 \mbox{ }\mathrm{for\ all}\mbox{ } s\in]t_{i-1},t_i]\}}]}{r_h}\\
&=&
\mathrm{O}\biggl(\frac{h \eta^2({{r_h}}/{4})}{r_h}\biggr)=
\mathrm{O}(h^{1-\alpha\beta
/2})
\end{eqnarray*}
and (\ref{PDXdmagSoglia}) is verified.
\end{pf}

\begin{pf*}{Proof of Proposition \protect\ref{teoStimIntQuart}}
\begin{eqnarray*}
&&\frac{\sum_i (\Delta_{i} X)^4I_{\{ (\Delta_{i} X)^2\leq r_h\}}
}{3h}\\
&&\quad = \frac{\sum_i (\Delta_{i} X_1)^4I_{\{ (\Delta_{i} X_1)^2\leq
4r_h\}}
}{3h}+
  \frac{1 }{3h}\sum_i (\Delta_{i} X_1)^4\bigl(I_{\{ (\Delta_{i}
X)^2\leq r_h\}} -I_{\{ (\Delta_{i} X_1)^2\leq4r_h\}}\bigr) \\
&&\qquad{}+
\sum_{k=1}^4 \pmatrix{
4 \cr k}
 \frac{\sum_i (\Delta_{i} X_1)^{4-k}(\Delta_{i}M)^k I_{\{
(\Delta_{i} X)^2\leq r_h\}} }{3h}\\
&&\quad :=
\sum_{j=1}^3 I_j(h).
\end{eqnarray*}
By Proposition 1 in \cite{mancini07}, $I_1(h)$ tends to $\int_0^T
\sigma_t^4 \,\mathrm{d}t$ in probability. We show here that the other terms
tend to zero in probability. Let us consider $I_2(h):= \frac{1
}{3h}\sum_i (\Delta_{i} X_1)^4(I_{\{ (\Delta_{i} X)^2\leq r_h\}}
-I_{\{ (\Delta_{i} X_1)^2\leq4r_h\}}) $: on $\{(\Delta_{i} X)^2\leq r_h,
(\Delta_{i} X_1)^2> 4r_h\},$ we have
\renewcommand{\theequation}{\arabic{equation}}
\setcounter{equation}{25}
\begin{equation} \label{DXpicceDXugr} \sqrt{r_h}\geq|\Delta_{i} X| >
|\Delta_{i} X_1|- |\Delta_{i}M| > 2\sqrt{r_h}- |\Delta_{i}M|,
\end{equation}
so $ |\Delta_{i}M|> \sqrt{r_h}$. Moreover,
if $|\Delta_{i} X_1|> 2\sqrt{r_h}$, then we necessarily have $\Delta
_{i} N\neq0$ since
\begin{equation} \label{DXugralloraDNpos} |\Delta_{i} X_0|+ |\Delta
_{i} J|\geq|\Delta_{i} X_1|> 2\sqrt{r_h}
\end{equation}
and, by (\ref{DXzeromagr}), a.s. for sufficiently small $h$, for all
$i=1,\ldots,n$, $|\Delta_{i} X_0|\leq\sqrt{r_h},$ thus $ |\Delta
_{i} J|> 2\sqrt{r_h}- |\Delta_{i} X_0| \geq
\sqrt{r_h}$. It follows that
\[
P\biggl\{\frac{1}{h}\sum_i (\Delta_{i} X_1)^4I_{\{(\Delta_{i}
X)^2\leq r_h, (\Delta_{i} X_1)^2> 4r_h\}}\neq0
\biggr\}
\leq n P\bigl\{|\Delta_{i}M|> \sqrt{r_h}, \Delta_{i} N\neq0\bigr\}
\rightarrow0,
\]
by Lemma
\ref{PDX0grandeePDtildeJ2grande0}.
On the other hand, for all $i=1,\ldots,n$ on 
$\{(\Delta_{i} X_1)^2\leq4r_h\},$ we have, for sufficiently small $h$,
$\Delta_{i} N=0$ because
\begin{equation} \label{DXuzeroalloraDNzero} |\Delta_{i} J|-|\Delta
_{i} X_0|\leq|\Delta_{i} X_1|\leq2\sqrt{r_h},
\end{equation}
so if $\Delta_{i} N\neq0,$ then a.s. for small $h$, we in fact have
$\Delta_{i} N=1$ and
$\Delta J_{ s} \geq1,$ by the definition of $J$.
Therefore, if $\Delta_{i} N\neq0,$ we would have $1\leq|\Delta_{i}
J|\leq2\sqrt{r_h}
+\sqrt{r_h}= 3\sqrt{r_h},$
which is impossible for small $h$. It follows that
\begin{eqnarray*}
\{(\Delta_{i} X)^2> r_h, (\Delta_{i} X_1)^2\leq4r_h\}&\subset&\{(
\Delta_{i} X_0+\Delta_{i}M)^2> r_h\}\\
&\subset&\biggl\{(
\Delta_{i} X_0)^2> \frac{r_h}{4}\biggr\}\cup
\biggl\{(\Delta_{i}M)^2> \frac{r_h}{4}\biggr\}.
\end{eqnarray*}
This implies, by (\ref{DXzeromagr}) and
(\ref{PDXdmagSoglia}), that a.s. as $h\rightarrow0,$
\begin{eqnarray*}
\frac{1}{h}\sum_i (\Delta_{i} X_1)^4I_{\{(\Delta_{i} X)^2> r_h,
(\Delta_{i} X_1)^2\leq4 r_h\}}&\leq&
\frac{\sum_i (\Delta_{i} X_0)^4 I_{\{(\Delta_{i}M)^2> r_h/4\}}}{h}
\\
&\leq&\Lambda^4 h\ln^2\frac1 h \sum_i I_{\{(\Delta_{i}M)^2> r_h/4\}}
\stackrel{P}{\rightarrow}0.
\end{eqnarray*}
%
We can conclude that $I_2(h)\stackrel{P}{\rightarrow}0$ as
$h\rightarrow0$. Now,
consider $I_3(h):= \sum_{k=1}^4 \left(
 {4 \atop k}
\right)I_{3, k}(h),$ where
\[
I_{3, k}(h):=\frac{1}{3h}\sum_i (\Delta_{i} X_1)^{4-k}(\Delta
_{i}M)^k I_{\{ (\Delta_{i} X)^2\leq r_h\}} ,\qquad
k=1,\ldots,4,
\]
is decomposable as
\begin{eqnarray}\label{sumDXukDJdj}&&\frac{1}{3h}\sum_i
(\Delta_{i} X_1)^{4-k}(\Delta_{i}M)^k I_{\{ (\Delta_{i} X)^2\leq
r_h, (\Delta_{i}M)^2\leq4r_h\}}
\nonumber
\\[-8pt]
\\[-8pt]
\nonumber
&&\quad{}+\frac{1}{3h}\sum_i (\Delta_{i} X_1)^{4-k}(\Delta_{i}M)^k I_{\{
(\Delta_{i} X)^2\leq r_h, (\Delta_{i}M)^2>
4r_h\}}.
\end{eqnarray}
We have, a.s. for small $h$, that for all $i$ on $\{ (\Delta_{i} X)^2
\leq
r_h, (\Delta_{i}M)^2> 4r_h\}$, $\Delta_{i} N\neq0$ since
\[
2\sqrt{r_h}-|\Delta_{i} X_1|< |\Delta_{i}M|- |\Delta_{i} X_1| 
\leq|\Delta_{i} X| \leq\sqrt{r_h}
\]
and then $|\Delta_{i} X_1| > \sqrt{r_h}$ and, similarly as
in (\ref{DXugralloraDNpos}), $|\Delta_{i} J|> 3\sqrt{r_h}/4$. So,
the probability
that the second term of (\ref{sumDXukDJdj}) differs from zero is
bounded by (\ref{DNneq0eDXdgeqr}) and tends to zero. As for the
first term, a.s. for sufficiently small $h$, for all $i$ on $\{
(\Delta_{i} X)^2\leq r_h, (\Delta_{i}M)^2\leq4r_h\},$ we have
$\Delta_{i} N=0$ because
\[
|\Delta_{i} X_1| - |\Delta_{i}M| \leq|\Delta_{i} X| \leq\sqrt{r_h},
\]
thus $|\Delta_{i} X_1|<3\sqrt{r_h}$ and we proceed as in (\ref
{DXuzeroalloraDNzero}).
So, 
the first term in (\ref{sumDXukDJdj}) is a.s. dominated by
\[
\frac{\sum_i |\Delta_{i} X_1|^{4-k}|\Delta_{i}M|^k I_{\{ \Delta
_{i} N=0, (\Delta_{i}M)^2\leq4r_h\}
}}{3h}\leq
\frac{\sum_i |\Delta_{i} X_0|^{4-k}|\Delta_{i}M|^k I_{\{ (\Delta
_{i}M)^2\leq4r_h\}}}{3h}.
\]
Now, for $k=4,$ 
we apply to $M$ property (C.19) in \cite{AitJacod08FavsIA}, Lemma 5,
with $\beta$ there being $\alpha$ here, $u_n= \sqrt{r_h}= h^{
\beta/2}$, $p=4$, $v_h= h^{\phi}$ for a proper exponent $\phi$ we specify
below and $\beta'=0$. Result (C.19) of \cite{AitJacod08FavsIA} then implies
that
\[
\frac1 h E\Biggl[\Biggl| \sum_{i=1}^n(\Delta_{i}M)^4I_{\{|\Delta
_{i}M|\leq2 \sqrt{r_h}\}} -
\sum_{v\leq T}|\Delta M_{v}|^4 I_{\{|\Delta M_{v}|\leq2
\sqrt{r_h}\}}\Biggr|\Biggr] \leq c h^{(\beta/2) (4-\alpha)-1}
\cdot\eta_{4,n},
\]
where $\eta_{4,n}= h(h^{\beta/2}
v_h)^{-\alpha} + h^2 h^{{\alpha\beta}/{2}}(h^{\beta/2}
v_h)^{-3\alpha} +h h^{{\alpha\beta}/{2}}(h^{\beta/2}
)^{-2\alpha} + (2 h^{\beta/2} )^{\alpha}  +\break h^{1/ 4}
h^{-((4-\alpha)/{4})\beta/2}+ v_h^{{(4-\alpha)}/{4}}$. As
soon as $\beta>1/(2-\alpha/2)$ and we choose $\phi\in\,]0,
\frac{1-\beta}{3}[$, so that for all $\alpha\in\,]0,2[$ we have $\phi
<(2/\alpha-\beta)/3$, it is guaranteed both that $h^{(\beta/2)
(4-\alpha)-1} \rightarrow0$ and that
$h^{({\beta}/{2})(4-\alpha)-1}\cdot\eta_{4,n} \rightarrow0$. Thus,
\[
\lim_h\frac{\sum_i |\Delta_{i}M|^4 I_{\{ (\Delta_{i}M)^2\leq4r_h\}}}{3h}=
\lim_h\frac{\sum_i \int_{t_{i-1}}^{t_i}\!\int_{|x|\leq2\sqrt
{r_h}}|x|^4 \mu(\mathrm{d}x,
\mathrm{d}t)}{3h}
\]
and
\begin{eqnarray*}
E\biggl[\sum_i \int_{t_{i-1}}^{t_i}\!\int_{|x|\leq2\sqrt{r_h}}|x|^4
\mu(\mathrm{d}x,
\mathrm{d}t)/3h\biggr]&=&\mathrm{O}\biggl(\int_{|x|\leq2\sqrt{r_h}}|x|^4 \nu(\mathrm{d}x)/h\biggr)\\
&=&
\mathrm{O} \bigl( h^{(\beta/2) (4-\alpha)-1}\bigr)\rightarrow0,
\end{eqnarray*}
given that $\beta >1/(2-\alpha /2)$.

To show, further, that the terms
\[
\frac{\sum_i |\Delta_{i}
X_0|^{4-k}|\Delta_{i}M|^k
I_{\{ (\Delta_{i}M)^2\leq4r_h\}}}{3h}
\]
tend to zero in probability for
$k=1,2,3,$ we use the fact that, by (\ref{DefDecompDXd}), each term is
dominated by (recall the notation in (\ref{defJh}))
\[
c   \frac{\sum_i |\Delta_{i} X_0|^{4-k}|\Delta_{i} M^{(h)}|^k }{3h}+
c  \frac{\sum_i |\Delta_{i} X_0|^{4-k}|h d(2\sqrt[4]{ r_h})|^k }{3h}.
\]
Now, a.s.
\begin{eqnarray*}
&&\frac{\sum_i |\Delta_{i} X_0|^{4-k}|h d(2\sqrt[4]{ r_h})|^k
}{3h}\\
&&\quad\leq
\biggl(h\ln\frac1 h \biggr)^{{(4-k)}/{2}} n h^{k-1} \biggl[\bigl|c+
r_h^{
{(1-\alpha)}/{4}}\bigr|^kI_{\{\alpha\neq1 \}}+ \ln^k\frac{1}{r_h^{1/4}}
I_{\{
\alpha=1\}}\biggr]\\
&&\quad\leq
c h^{k/2} \biggl(\ln\frac1 h \biggr)^{{(4-k)}/{2}} +
c h^{k/2} \biggl(\ln\frac1 h\biggr)^{{(4-k)}/{2}}  r_h^{k{(1-\alpha)}/{4}}
+ h^{ h/ 2} \ln^{2-k/2}\frac{1}{hr_h^{1/4}}\\
&&\quad= \mathrm{o}(1) +c h^{k[1/ 2 +\beta{(1-\alpha)}/{4}]} \log^{{(4-k)}/{2}}\frac1 h \rightarrow0
\end{eqnarray*}
for all $k=1,2,3$ as
$r_h= h^\beta$, $\beta\in\,]0,1[$. 
As for
\begin{equation} \label{sunDXzDXdm} \frac{\sum_i |\Delta_{i}
X_0|^{4-k}|\Delta_{i} M^{(h)}|^k
}{3h},
\end{equation}
we need to deal separately with each of $k=1,2,3$. Note that since
$a$ and $\sigma$ are locally bounded on $\Omega\times[0,T]$, we
can assume that they are bounded without loss of generality,
so $E[(\int_{t_{i-1}}^{t_i}\sigma_s\,\mathrm{d}W_s)^{2k}]=\mathrm{O}(h^k)$ for each $k=1,2,3$,
using, for instance, the Burkholder inequality~\cite{Protter05}, page
226,
and a.s. 
$(\int_{t_{i-1}}^{t_i}a_s \,\mathrm{d}s)^{2k}
=\mathrm{o}(h^k)$. Therefore, $E[(\Delta_{i} X_0)^{2k}]= \mathrm{O}(h^k)$ for each of $k=1,2,3$.
For $k=1,$ the expected value of (\ref{sunDXzDXdm}) is bounded by
$({n}/{3h}) \sqrt{E[(\Delta_{i} X_0)^6]}\times\sqrt{E(\Delta_{i}
M^{(h)})^2} = \mathrm{O}(r_h^{(1 /4)(1-\alpha/2)})$ and thus
tends to zero as $h\rightarrow0$.
As for $k=2,$
\begin{equation} \label{sumDXzqDXdqsuh}
\frac{\sum_i (\Delta_{i} X_0)^{2}(\Delta_{i} M^{(h)})^2 }{h}\leq h
\ln\frac1 h \frac{\sum_i
(\Delta_{i} M^{(h)})^{2}}{h},
\end{equation}
whose expected value is given by
\[
\ln\frac1 h \eta^2(2r_h^{1/ 4})\rightarrow0
\]
as $h\rightarrow0$ since $r_h=h^\beta$, with $\beta>0$. Concerning
$k=3,$ we have
\[
\frac{\sum_i |\Delta_{i} X_0| |\Delta_{i} M^{(h)}|^3 }{h}\leq\frac
c h \sum_i (\Delta_{i} X_0)^{2}\bigl(\Delta_{i} M^{(h)}
\bigr)^2 + \frac c h \sum_i \bigl(\Delta_{i} M^{(h)}\bigr)^4,
\]
so that this step is reduced to the steps with $k=2,4$ which we dealt
with previously.
\end{pf*}

\begin{pf*}{Proof of Theorem \protect\ref{teoSpeedCvSumDXdq}} 
Let us define $K_{ni}:= (\int_{t_{i-1}}^{t_i}\!\int_{|x|\leq
\varepsilon} x \tilde\mu
(\mathrm{d}x,\mathrm{d}t)-h\int_{\varepsilon<|x|\leq1} x \nu(\mathrm{d}x) )^2 $. We
apply the
Lindeberg--Feller theorem to the double array sequence $H_{n i}$
given by the normalized versions of the variables $K_{ni}$,
$i=1,\ldots,n$, and $n=T/h$. Using relations (\ref{OrdGrandIntPotxdnu}),
we have
\renewcommand{\theequation}{\arabic{equation}}
\setcounter{equation}{31}
\begin{eqnarray} \label{Enj}
E[K_{ni}]&=&h\ell_{2,h}\varepsilon^{2-\alpha}+ \biggl(h
\int_{
\varepsilon<|x|\leq1} x \nu(\mathrm{d}x) \biggr)^2
\nonumber
\\[-8pt]
\\[-8pt]
\nonumber
& =& h\ell
_{2,h}\varepsilon^{2-\alpha} +
\ell_{1,h}^2 h^2 \biggl[(c+\varepsilon^{1-\alpha})^2 I_{\{\alpha\neq
1\}} +
\biggl(\ln^2  \frac{1}{\varepsilon}\biggr)I_{\{\alpha= 1\}}\biggr].
\end{eqnarray}
%
Taking $\varepsilon=h^u$ and any $u\in\,]0, 1/2],$ we obtain that
\begin{eqnarray*}
v^2_{ni}:= \operatorname{var}[K_{ni}]&=& E\biggl[\biggl(\int
_{t_{i-1}}^{t_i}\!\int_{|x|\leq\varepsilon
} x \tilde\mu(\mathrm{d}x,\mathrm{d}t)-
h\int_{\varepsilon<|x|\leq1} x \nu(\mathrm{d}x) \biggr)^4\biggr]\\
&&{}- E^2_{ni}
\sim
h\int_{|x|\leq\varepsilon} x^4 \nu(\mathrm{d}x) = h \ell_{4,h} \varepsilon
^{4-\alpha}
\end{eqnarray*}
as
$h\rightarrow0$. Then, consider
\[
H_{n i}:= \frac{K_{ni} -E[K_{ni}]}{\sqrt n v_{ni}}
\sim\frac{K_{ni} - h\ell_{2,h}\varepsilon^{2-\alpha}
- \ell_{1,h}^2 h^2 [(c+\varepsilon^{1-\alpha})^2I_{\{\alpha
\neq1\}}+
(\ln^2{1}/{\varepsilon}) I_{\{\alpha=1\}}]}{\sqrt T
\sqrt{\ell
_{4,h}} \varepsilon^{2-\alpha/2}}.
\]
%
We now show that for any $\delta>0$, there exists a $q>1$ such that
\begin{equation} \label{Lind}
\sum_{i=1}^n E\bigl[H_{n i}^2 I_{\{|H_{n i}|>\delta\}}\bigr]\leq c\varepsilon
^{
{\alpha}/{(2q)}}\rightarrow0
\end{equation}
as $h\rightarrow0$, so the Lindeberg condition
is satisfied and implies that
\begin{equation} \label{TeoFeller}
\sum_{i=1}^n H_{n i} \stackrel{d}\rightarrow N(0,1).
\end{equation}
Noting that $h/\varepsilon^{2-\alpha/2}$ and $(h\varepsilon
^{1-\alpha})/(\varepsilon^{2-\alpha
/2}) I_{\{\alpha\neq1\}} + (h\ln^2(1/\varepsilon))/(\varepsilon
^{2-\alpha/2}) I_{\{
\alpha=1\}}$ tend to zero as $h\rightarrow0$, (\ref{TeoFeller})
leads to (\ref
{SpeedCvSumDXdq}).
To show inequality (\ref{Lind}), consider
\begin{equation} \label{HolderperLind}
n E\bigl[H_{n 1}^2 I_{\{|H_{n 1}|>\delta\}}\bigr]\leq n E^{1/ p}[H_{n
1}^{2p}]P^{1 /q}\{|H_{n 1}|>\delta\},
\end{equation}
as for the last factor above, we note that $|H_{n 1}|>\delta$ if and
only if either
\begin{eqnarray*}
K_{n1}&<&h\ell_{2,h}\varepsilon^{2-\alpha} + \ell_{1,h}^2 h^2
\biggl[(c+\varepsilon
^{1-\alpha})^2 I_{\{\alpha\neq1\}} + \biggl(\ln^2\frac
{1}{\varepsilon}\biggr) I_{\{
\alpha=1\}}\biggr]
-\delta\sqrt{T \ell_{4,h}}\varepsilon^{2-\alpha/2}\\
&=&
\varepsilon^{2-\alpha/2}
\bigl(\mathrm{o}(1)-c \delta\bigr),
\end{eqnarray*}
where $c$ denotes a generic constant, or
\[
K_{n1} >
h\ell_{2,h}\varepsilon^{2-\alpha}+\ell_{1,h}^2 h^2
\biggl[(c+\varepsilon^{1-\alpha})^2
I_{\{\alpha\neq1\}} + \biggl(\ln^2\frac{1}{\varepsilon}\biggr)
I_{\{\alpha=1\}}\biggr]+c \delta\varepsilon^{2-\alpha/
2}=\mathrm{O}(\varepsilon^{2-
\alpha/2 }).
\]
However, $K_{n1}\geq0$, while for sufficiently small
$h$, the right-hand term of the first inequality above is strictly
negative, therefore $|H_{n 1}|>\delta$ if and only if $K_{n1}>c
\varepsilon
^{2-
\alpha/2 }$, that is, either
\[
-c \varepsilon^{1-\alpha/4}\sim h(c+\varepsilon^{1- \alpha}
)I_{\{\alpha\neq1\}}
+ I_{\{\alpha=1\}} h\ln\frac1 \varepsilon- c\varepsilon^{1-\alpha/4}>
\int_0^{t_1}\!\!\! \int_{|x|\leq\varepsilon} x \tilde\mu(\mathrm{d}x,\mathrm{d}t)
\]
or, for sufficiently small $h$,
$
\int_0^{t_1}\!\!\! \int_{|x|\leq\varepsilon} x \tilde\mu(\mathrm{d}x,\mathrm{d}t) >
c\varepsilon^{1-
\alpha/4},$ and so $|H_{n 1}|>\delta$ if and only if
\[
\biggl|\int_0^{t_1}\!\!\! \int_{|x|\leq\varepsilon} x \tilde\mu
(\mathrm{d}x,\mathrm{d}t)\biggr|> c \varepsilon
^{1-\alpha/4}.
\]
This entails that for sufficiently small $h$,
\begin{eqnarray*}
P\{|H_{n 1}|>\delta\}&=&
P\biggl\{\biggl|\int_0^{t_1}\!\!\! \int_{|x|\leq\varepsilon} x \tilde\mu
(\mathrm{d}x,\mathrm{d}t)\biggr|>
c \varepsilon^{1-\alpha/4}\biggr\}\\
&\leq&
c \frac{E[|\int_0^{t_1}\!\!\! \int_{|x|\leq\varepsilon} x \tilde\mu
(\mathrm{d}x,\mathrm{d}t)|^2]}{\varepsilon
^{2-\alpha/2}}= h^{1-{(\alpha u)}/{2}}\rightarrow0.
\end{eqnarray*}
The first two factors of the right-hand side of (\ref{HolderperLind})
are dominated by
\begin{eqnarray*}
&&c n \frac{E^{1 /p}[(K_{n1} - h \ell_{2,h}\varepsilon
^{2-\alpha}
- h^2 \ell_{1,h}^2 [(c+\varepsilon^{1-\alpha})^2 I_{\{\alpha
\neq1\}} + (\ln
^2{1}/{\varepsilon}) I_{\{\alpha=1\}}])^{2p}
]}{\varepsilon^{4-\alpha}}
\\
&&\quad\leq
c n \frac{E^{1/ p}[K_{n1}^{2p}] + (h\varepsilon
^{2-\alpha})^2 +
h^4 (1-\varepsilon^{1-\alpha})^4 +
h^4 \ln^4{1}/{\varepsilon}}{\varepsilon^{4-\alpha}}.
\end{eqnarray*}
The last three terms give no contribution to (\ref{HolderperLind})
since
\[
n \frac{(h\varepsilon^{2-\alpha})^2 + h^4 (1-\varepsilon^{1-\alpha
})^4+h^4 \ln^4
{1}/{\varepsilon}}{\varepsilon^{4-\alpha}} h^{(1-{\alpha
u}/{2})(1 /q)} \rightarrow0.
\]
On the other hand, by choosing, for example, $p=5/4$, we have 
%
\[
E[K_{n1}^{2p}] = \mathrm{O}(h\varepsilon^{5-\alpha}),
\]
so we are left to deal with $n \frac{(h\varepsilon^{5-\alpha
})^{1/ p} }{\varepsilon
^{4-\alpha}} h^{(1-{\alpha u}/{2})(1 /q)} = \varepsilon
^{{\alpha}/{(2q)}}$ and the inequality in (\ref{Lind}) is proved.
\end{pf*}

\begin{Lemma} \label{IndDxqleqrvsIndDxdqleqr}
As $h\rightarrow0$, if $r_h\rightarrow0$, $n=T/h$ and
$\sup_{i=1,\ldots,n} |a_{hi}| = \mathrm{O}(r_h),$
then
\[
\sum_i |a_{hi}|I_{\{ (\Delta_{i} X)^2\leq r_h\}} - \sum_i |a_{hi}|
I_{\{(\Delta_{i}M)^2\leq4 r_h, \Delta_{i} N=0\}
}\stackrel{P}{\rightarrow}0.
\]
\end{Lemma}

\begin{pf}
On $\{(\Delta_{i} X)^2\leq r_h\},$ we have $ |\Delta_i L |-
|\Delta_{i} X_0| \leq|\Delta_{i} X|\leq\sqrt{r_h}$ and, thus, by
(\ref{DXzDivhlnunosuh}),
for small~$h$, $ |\Delta_i L |\leq2\sqrt{r_h}$, so that a.s.
\[
\lim _{h\rightarrow0} \sum_i |a_{hi}|I_{\{ (\Delta_{i}
X)^2\leq r_h\}} \leq\lim _{h\rightarrow0}
\sum_i |a_{hi}| I_{\{(\Delta_i L)^2 \leq4 r_h\}}.
\]
However,
\renewcommand{\theequation}{\arabic{equation}}
\setcounter{equation}{35}
\begin{equation} \label{sumaniIDNneqzerotendszero}
\sum_i |a_{hi}| I_{\{(\Delta_i L)^2 \leq4 r_h,
\Delta_{i} N\neq0\}}\leq\sup_i |a_{hi}| N_T \stackrel
{\mathit{a.s.}}\rightarrow0
\end{equation}
as $h\rightarrow0$ and thus a.s.
\[
\lim _{h\rightarrow0} \sum_i |a_{hi}|I_{\{ (\Delta_{i}
X)^2\leq r_h\}} \leq
\lim _{h\rightarrow0} \sum_i |a_{hi}| I_{\{(\Delta_i L)^2
\leq4 r_h, \Delta_{i} N
=0\}}=
\lim _{h\rightarrow0} \sum_i |a_{hi}| I_{\{(\Delta
_{i}M)^2\leq4 r_h, \Delta_{i} N=0\}}.
\]
We now show that, on the other hand, the positive quantity
\[
\lim _{h\rightarrow0} \sum_i|a_{hi}| \bigl(I_{\{(\Delta_i L)^2
\leq4 r_h, \Delta_{i} N
=0\}}-
I_{\{(\Delta_{i} X)^2\leq r_h\}} \bigr)= 0 \qquad  \mbox{a.s.}
\]
In fact,
\begin{eqnarray*}
&&\{(\Delta_i L)^2 \leq4 r_h, \Delta_{i} N=0\}-\{(\Delta_{i} X)^2\leq
r_h\}\\
&&\quad =
\{(\Delta_i L)^2 \leq4 r_h, \Delta_{i} N=0, (\Delta_{i} X)^2> r_h\}\\
&&\qquad\subset
\bigl\{|\Delta_i L| \leq2\sqrt{r_h}, \Delta_{i} N=0, |\Delta_{i}
X_0|+|\Delta_{i}M| > \sqrt{r_h}\bigr\}\\
&&\qquad\subset
\bigl\{|\Delta_{i} X_0| >   \sqrt{r_h}/2\bigr\} \cup
\bigl\{|\Delta_{i}M| \leq2\sqrt{r_h}, |\Delta_{i}M| >   \sqrt
{r_h}/2\bigr\}.
\end{eqnarray*}
Since, by (\ref{DXzeromagr}), a.s. for sufficiently small $h$
$ \sum_i|a_{hi}| I_{\{|\Delta_{i} X_0|> \sqrt{r_h}/2\}}=0,$
we a.s. have
\[
\lim _{h\rightarrow0} \sum_i|a_{hi}| \bigl(I_{\{(\Delta_i L)^2
\leq4 r_h,
\Delta_{i} N=0\}}- I_{\{(\Delta_{i} X)^2\leq r_h\}} \bigr)\leq
\lim _{h\rightarrow0} \sum_i|a_{hi}| I_{\{|\Delta_{i}M| \leq
2\sqrt{r_h},|\Delta_{i}M| >
\sqrt{r_h}/2\}} ;
\]
however, by Remark \ref{JumpsLeqrOnefourth}, as $h\rightarrow0,$ 
%
\begin{eqnarray*}
E\biggl[\sum_i|a_{hi}| I_{\{|\Delta_{i}M| \leq2\sqrt{r_h}, |\Delta
_{i}M|> \sqrt{r_h}/2\}}\biggr]&\leq&
\mathrm{O}(r_h)n P\bigl\{|\Delta_{i}M| \leq2\sqrt{r_h},|\Delta_{i}M|> \sqrt
{r_h}/2\bigr\}\\
&\leq&
\mathrm{O}(r_h)n P\bigl\{|\Delta_{i}M|I_{\{|\Delta_{i}M| \leq2\sqrt{r_h}\}}>
\sqrt{r_h}/2\bigr\}\\
&\leq&
\mathrm{O}(r_h)n \frac{E[(\Delta_{i}M)^2I_{\{|\Delta_{i}M| \leq2\sqrt{r_h}\}
}]}{r_h} \\
&=& \mathrm{O}(r_h) n
\frac{h \eta^2(2r_h^{c1 /4})}{r_h} \rightarrow0.
\end{eqnarray*}
\upqed\end{pf}

\begin{Lemma}\label{LemSumDXdqdoveDXqleqrAndDXdqleqr}
Under the assumptions of Theorem \ref{corVelCvIV}, for all $\alpha\in[0,2[,$
\renewcommand{\theequation}{\arabic{equation}}
\setcounter{equation}{36}
\begin{eqnarray} \label{I4}
\sum_{i=1}^n (\Delta_{i}M)^2I_{\{(\Delta_{i}M)^2\leq r_h/16\}}-
\mathrm{o}_P(h^{1-\alpha/2})&\leq&
\sum_{i=1}^n (\Delta_{i}M)^2I_{\{(\Delta_{i} X)^2\leq r_h, (\Delta
_{i}M)^2\leq4r_h\}}
\nonumber
\\[-8pt]
\\[-8pt]
\nonumber
&\leq&\sum_{i=1}^n (\Delta_{i}M)^2I_{\{(\Delta_{i}M)^2\leq9r_h/4\}}+
\mathrm{o}_P(h^{1-\alpha
/2}) \qquad a.s.
\end{eqnarray}
\end{Lemma}
\begin{pf}Let us first deal with $\sum_{i=1}^n (\Delta
_{i}M)^2I_{\{(\Delta_{i} X)^2>r_h, (\Delta_{i}M)^2
\leq4r_h\}}$.

As in (\ref{doveDXqgrr}), on $ \{(\Delta_{i} X)^2>r_h\},$
we have either $|\Delta_{i} J
|> \sqrt{r_h}/4 $ or $ |\Delta_{i}M|> \sqrt{r_h}/4 $, so
\begin{eqnarray*}
&&\sum_{i=1}^n (\Delta_{i}M)^2I_{\{(\Delta_{i} X)^2>r_h, (\Delta
_{i}M)^2\leq4r_h\}}\\
&&\quad\leq
\sum_{i=1}^n (\Delta_{i}M)^2I_{\{(\Delta_{i} X)^2>r_h,\Delta_i
J\neq0, (\Delta_{i}M)^2\leq4r_h\}}\\
&&\qquad{} +
\sum_{i=1}^n (\Delta_{i}M)^2I_{\{(\Delta_{i} X)^2>r_h,(\Delta
_{i}M)^2> {r_h}/{16}, (\Delta_{i}M)^2\leq4r_h\}
}.
\end{eqnarray*}
%
However,
\[
E\biggl[\frac{\sum_{i=1}^n (\Delta_{i}M)^2I_{\{(\Delta_{i}M)^2\leq
4r_h, \Delta_{i} N\neq0 \}
}}{h^{1-\alpha/2}}\biggr]
= \mathrm{O}\biggl(\frac{h\eta^2(r_h^{1/ 4})N_T}{h^{1-\alpha/2}}
\biggr)\rightarrow
0,
\]
so
\renewcommand{\theequation}{\arabic{equation}}
\setcounter{equation}{37}
\begin{eqnarray}\label{maggioraz}
&&\sum_{i=1}^n (\Delta_{i}M)^2I_{\{(\Delta_{i}
X)^2>r_h, (\Delta_{i}M)^2\leq
4r_h\}}\nonumber\\
&&\quad\leq
\mathrm{o}_P(h^{1-\alpha/2})+ \sum_{i=1}^n (\Delta_{i}M)^2I_{\{(\Delta
_{i}M)^2\leq4r_h, (\Delta_{i}M)^2>r_h
/16\}}
\nonumber
\\[-8pt]
\\[-8pt]
\nonumber
&&\quad=\mathrm{o}_P(h^{1-\alpha/2})+ \sum_{i=1}^n (\Delta_{i}M)^2I_{\{(\Delta
_{i}M)^2\leq4r_h\}}- \sum
_{i=1}^n (\Delta_{i}M)^2I_{\{(\Delta_{i}M)^2\leq4r_h, (\Delta
_{i}M)^2\leq r_h/16\}}
\\
&&\quad= \mathrm{o}_P(h^{1-\alpha/2})+ \sum_{i=1}^n (\Delta_{i}M)^2I_{\{(\Delta
_{i}M)^2\leq4r_h\}}- \sum
_{i=1}^n (\Delta_{i}M)^2I_{\{(\Delta_{i}M)^2\leq r_h/16\}}.\nonumber
\end{eqnarray}
Now, consider $ \sum_{i=1}^n (\Delta_{i}M)^2I_{\{(\Delta_{i}M)^2\leq
4r_h, (\Delta_{i}M)^2>9r_h/4\}
}$: on
$\{2\sqrt{r_h}\geq|\Delta_{i}M|> \frac3 2 \sqrt{r_h}\},$
either $\Delta_{i} N\neq0$, in which case
\[
\frac{\sum_{i=1}^n (\Delta_{i}M)^2I_{\{(\Delta_{i}M)^2\leq4r_h,
\Delta_{i} N\neq0\}}}{h^{1-\alpha
/2}}\stackrel{P}{\rightarrow} 0,
\]
as before, or $\Delta_{i} N=0$, in which case $ |\Delta_{i} X|>
|\Delta_{i}M|- |\Delta_{i} X_0|> \frac3 2
\sqrt{r_h}- \frac1 2 \sqrt{r_h}=
\sqrt{r_h}$, so
\[
\sum_{i=1}^n (\Delta_{i}M)^2I_{\{(\Delta_{i} X)^2>r_h, (\Delta
_{i}M)^2\leq4r_h\}}+
\mathrm{o}_P(h^{1-\alpha/2})\geq\sum_{i=1}^n (\Delta_{i}M)^2I_{\{(\Delta
_{i}M)^2\leq4r_h, (\Delta_{i}M)^2
>9r_h/4\}}.
\]
Therefore,
\begin{eqnarray} \label{minoraz}
&&\sum_{i=1}^n (\Delta_{i}M)^2I_{\{(\Delta_{i} X)^2>r_h, (\Delta
_{i}M)^2\leq4r_h\}}
\nonumber
\\[-8pt]
\\[-8pt]
\nonumber
&&\quad\geq-
\mathrm{o}_P(h^{1-\alpha/2})+ \sum_{i=1}^n (\Delta_{i}M)^2I_{\{(\Delta
_{i}M)^2\leq4r_h\}}- \sum
_{i=1}^n (\Delta_{i}M)^2I_{\{(\Delta_{i}M)^2\leq9r_h/4\}}.
\end{eqnarray}
Now combining (\ref{maggioraz}) and (\ref{minoraz}), we obtain (\ref
{I4}) since
\begin{eqnarray*}
&&\sum_{i=1}^n (\Delta_{i}M)^2I_{\{(\Delta_{i} X)^2\leq r_h, (\Delta
_{i}M)^2\leq4r_h\}}\\
&&\quad=
\sum_{i=1}^n (\Delta_{i}M)^2I_{\{(\Delta_{i}M)^2\leq4r_h\}}
- \sum_{i=1}^n (\Delta_{i}M)^2I_{\{(\Delta_{i} X)^2>r_h, (\Delta
_{i}M)^2\leq4r_h\}}.
\end{eqnarray*}
\upqed\end{pf}

\begin{pf*}{Proof of Theorem \protect\ref{corVelCvIV}} Note that under $\beta>
\frac{1}{2-\alpha/2}$, the assumptions of Proposition
\ref{teoStimIntQuart} are satisfied. Since $X=X_1+ M$, we
decompose
\renewcommand{\theequation}{\arabic{equation}}
\setcounter{equation}{39}
\begin{eqnarray}\label{CLTterm}
\frac{\hat{\mathit{IV}}_h - \mathit{IV}}{\sqrt{2h\hat{\mathit{IQ}}_h}}&=&\frac{\sum_{i=1}^n
(\Delta_{i} X)^2
I_{\{ (\Delta_{i} X)^2\leq r_h\}} -\mathit{IV}}{
\sqrt{2h}\sqrt{{\sum_i (\Delta_{i} X)^4I_{\{ (\Delta_{i}
X)^2\leq r_h\}} }/{3h}}}
\nonumber
\\[-8pt]
\\[-8pt]
\nonumber
&=&
\frac{\sum_{i=1}^n (\Delta_{i} X_1)^2I_{\{(\Delta_{i} X_1)^2\leq
4r_h\}}-\mathit{IV}}{
\sqrt{({2}/{3})\sum_i (\Delta_{i} X)^4I_{\{ (\Delta_{i} X)^2\leq
r_h\}} }}\\
\label{dimCor4term} &&{}
+\frac{\sqrt{2hIQ}}{\sqrt{({2}/{3})\sum_i (\Delta_{i} X)^4I_{\{
(\Delta_{i} X)^2\leq r_h\}} }}\nonumber\\
&&\quad{}\times\biggl[ \frac{\sum_{i=1}^n (\Delta_{i} X_1)^2( I_{\{(\Delta_{i}
X)^2\leq r_h\}} - I_{\{
(\Delta_{i} X_1)^2\leq4r_h\}})}{
\sqrt{2h \mathit{IQ}}}
\nonumber
\\[-8pt]
\\[-8pt]
\nonumber
&&\hspace*{18pt}{}+
2\frac{\sum_{i=1}^n \Delta_{i} X_1\Delta_{i}MI_{\{ (\Delta_{i}
X)^2\leq r_h\}} }{\sqrt{2h \mathit{IQ}}}
+ \frac{\sum_{i=1}^n (\Delta_{i}M)^2I_{\{ (\Delta_{i} X)^2\leq r_h\}
} }{\sqrt{2h \mathit{IQ}}}\biggr]\qquad\\
&:=& \sum
_{j=1}^4 I_j(h).\nonumber
\end{eqnarray}
The proof of \cite{mancini07}, Theorem 2, shows that
$I_1(h)$ converges stably in law to a standard Gaussian random variable.
To show that the remaining terms either tend to zero or to infinity,
we can assume without loss of generality that both $a $ and
$\sigma$ are bounded a.s. If $(\Delta_{i} X)^2\leq r_h$ and $(\Delta
_{i} X_1)^2
> 4r_h$, then $ |\Delta_{i}M|> \sqrt{r_h}$ and $\Delta_{i} N\neq0$,
exactly as for
$I_2(h)$ in Proposition \ref{teoStimIntQuart}. It follows that
\begin{eqnarray*}
&&P\biggl\{ \frac{\sum_{i=1}^n (\Delta_{i} X_1)^2I_{\{(\Delta_{i}
X)^2\leq r_h, (\Delta_{i} X_1)^2> 4r_h\}}
}{\sqrt{2h \mathit{IQ}}} \neq0 \biggr\}\\
&&\quad\leq
nP\bigl\{\Delta_{i} N\neq0,|\Delta_{i}M|> \sqrt{r_h}\bigr\} \rightarrow0,
\end{eqnarray*}
by (\ref{DNneq0eDXdgeqr}).
The main factor of the remaining part of $I_2(h)$ is
\[
\frac{\sum_{i=1}^n (\Delta_{i} X_1)^2I_{\{(\Delta_{i} X)^2> r_h,
(\Delta_{i} X_1)^2\leq4r_h\}} }{\sqrt
{2h \mathit{IQ}}}.
\]
We recall that on $\{|\Delta_{i} X_1|\leq2\sqrt{r_h}\},$ we have
$\Delta_{i} N=0$, thus $(\Delta_{i} X_1)^2
=(\Delta_{i} X_0)^2$.
Moreover,
\[
\frac{\sum_{i=1}^n (\int_{t_{i-1}}^{t_i}a_u \,\mathrm{d}u)^2
I_{\{(\Delta_{i} X)^2> r_h, (\Delta_{i} X_1)^2\leq4r_h\}} }{\sqrt
{2h \mathit{IQ}}}= \mathrm{O}_P\bigl(\sqrt h\bigr) \rightarrow0
\]
and, by (\ref{SumIDXqmagr}),
\begin{eqnarray*}
\frac{1}{\sqrt{2h \mathit{IQ}}}\sum_{i=1}^n \int_{t_{i-1}}^{t_i}a_u \,\mathrm{d}u \int
_{t_{i-1}}^{t_i}\sigma_u \,\mathrm{d}W_u\,
I_{\{(\Delta_{i} X)^2> r_h, (\Delta_{i} X_1)^2\leq4r_h\}} &\leq& c
\sqrt h \sqrt{h\ln\frac1 h}\sum_{i=1}^n
I_{\{(\Delta_{i} X)^2> r_h\}}
\\
&=& \mathrm{O}\Biggl(h^{1-\alpha\beta/2}\sqrt{\ln\frac1 h}\Biggr) \rightarrow0.
\end{eqnarray*}
Therefore, in probability,
\[
\lim_{h\rightarrow0} I_2(h)=
\lim_{h\rightarrow0}-
\frac{\sum_{i=1}^n (\int_{t_{i-1}}^{t_i}\sigma_u \,\mathrm{d}W_u)^2
I_{\{(\Delta_{i} X)^2>r_h,(\Delta_{i} X_1)^2
\leq4r_h\}} }{\sqrt{2h \mathit{IQ}}}.
\]

We now show that term $I_3(h)/2$ in (\ref{dimCor4term}) tends to zero
in probability.
First, recall that $\Delta_{i} X_1= \Delta_{i} X_0+ \Delta_{i} J$
and, within the sum
$ \sum_{i=1}^n \Delta_{i} J\Delta_{i}MI_{\{(\Delta_{i} X)^2\leq
r_h\}}/\sqrt h$, the
term $i$ contributes only when $\Delta_{i} N\neq0$, in which case we
also have
$(\Delta_{i} X_1)^2> 4r_h$ and thus
$|\Delta_{i}M|> \sqrt{r_h},$ as in (\ref{DXpicceDXugr}). That implies
\[
P\biggl\{\frac{\sum_{i=1}^n \Delta_{i} J\Delta_{i}MI_{\{ (\Delta
_{i} X)^2\leq r_h\}} }{\sqrt{2h \mathit{IQ}}}\neq
0\biggr\}\leq
n P\bigl\{\Delta_{i} N\neq0 , |\Delta_{i}M|> \sqrt{r_h}\bigr\} \rightarrow0.
\]
As for $\frac{\sum_{i=1}^n
\Delta_{i} X_0\Delta_{i}MI_{\{ (\Delta_{i} X)^2\leq r_h\}} }{\sqrt
{h }}$, as in the proof of Lemma
\ref{IndDxqleqrvsIndDxdqleqr}, we have
\begin{equation} \label{DoveDXqleqrhoDXdqleq4r}
\frac{\sum_{i=1}^n \Delta_{i} X_0\Delta_{i}MI_{\{ (\Delta_{i}
X)^2\leq r_h\}} }{\sqrt{h }}=\frac{\sum_{i=1}^n
\Delta_{i} X_0\Delta_{i}M
I_{\{(\Delta_{i} X)^2\leq r_h, (\Delta_{i} L)^2\leq4r_h\}}}{\sqrt{h }}.
\end{equation}
However, since both
$P\{\frac{1}{\sqrt h} \sum_{i=1}^n
\Delta_{i} X_0\Delta_{i}MI_{\{(\Delta_{i} X)^2\leq r_h, (\Delta_{i}
L)^2\leq4r_h, \Delta_{i} N\neq0\}}\neq0 \}$
and
$P\{\frac{1}{\sqrt h}\times \sum_{i=1}^n \Delta_{i} X_0\Delta_{i}M
I_{\{(\Delta_{i} X)^2\leq r_h, (\Delta_{i}M)^2\leq4r_h, \Delta_{i}
N\neq0\}}\neq0 \}$
are dominated by\break $ n P\{\Delta_{i} N\neq0,  (\Delta_{i}M)^2> cr_h\}
\rightarrow0,$ we have
%
\begin{eqnarray*}
&&\lim_h \frac{1}{\sqrt h} \sum_{i=1}^n \Delta_{i} X_0\Delta_{i}M
I_{\{(\Delta_{i} X)^2\leq r_h, (\Delta_{i} L)^2\leq4r_h\}}\\
&&\quad=
\lim_h \frac{1}{\sqrt h} \sum_{i=1}^n \Delta_{i} X_0\Delta
_{i}MI_{\{(\Delta_{i} X)^2\leq r_h, (\Delta_{i} L)^2
\leq4r_h, \Delta_{i} N=0\}}
\\
&&\quad= \lim_h \frac{1}{\sqrt h} \sum_{i=1}^n \Delta_{i} X_0\Delta
_{i}MI_{\{(\Delta_{i} X)^2\leq r_h,
(\Delta_{i}M)^2\leq4r_h, \Delta_{i} N=0\}}\\
&&\quad=
\lim_h \frac{1}{\sqrt h} \sum_{i=1}^n \Delta_{i} X_0\Delta
_{i}MI_{\{(\Delta_{i} X)^2\leq r_h,
(\Delta_{i}M)^2\leq4r_h\}}.
\end{eqnarray*}
Moreover,
by the Cauchy--Schwarz inequality, we have
\begin{eqnarray}\label{QVofSmallJs}
&&\frac{\sum_{i=1}^n \int_{t_{i-1}}^{t_i}a_u \,\mathrm{d}u\, \Delta_{i}MI_{\{
(\Delta_{i} X)^2\leq r_h, (\Delta_{i}M)^2\leq4r_h
\}}}{\sqrt{h }}\nonumber\\
&&\quad\leq
\frac{\sqrt{\sum_{i=1}^n (\int_{t_{i-1}}^{t_i}a_u \,\mathrm{d}u)^2 }}{\sqrt
{h}} \sqrt{\sum
_{i=1}^n (\Delta_{i}M)^2
I_{\{ (\Delta_{i}M)^2\leq4 r_h\}}}
\\
&&\quad\leq c 
 \sqrt{\sum_{i=1}^n (\Delta_{i}M)^2I_{\{ (\Delta_{i}M)^2\leq4
r_h\}}},\nonumber
\end{eqnarray}
which tends to zero in
probability since, by Remark \ref{JumpsLeqrOnefourth}, as
$h\rightarrow0,$
\begin{equation} \label{VQThresholdedJumps}E\Biggl[\sum_{i=1}^n
(\Delta_{i}M)^2I_{\{ (\Delta_{i}M)^2\leq4 r_h\}}\Biggr]
=\int_0^T\!\!\! \int_{|x| \leq2 r_h^{1/4}} x^2 \nu(\mathrm{d}x)= T\eta^2(r_h
^{1/4})\rightarrow0.
\end{equation}
%
On the other hand,
\begin{eqnarray}\label{sumsigmaWDXdIDXdpiccolodivisosqrth}
&&\frac{1}{\sqrt h} \sum_{i=1}^n \biggl(\int_{t_{i-1}}^{t_i}
\sigma_u \,\mathrm{d}W_u
\biggr)\Delta_{i}MI_{\{(\Delta_{i} X)^2\leq r_h, (\Delta_{i}M)^2\leq4r_h\}
}  \nonumber\\
&&\quad=
\frac{1}{\sqrt h} \sum_{i=1}^n \biggl(\int_{t_{i-1}}^{t_i}
\sigma_u \,\mathrm{d}W_u
\biggr)\Delta_{i} M^{(h)}I_{\{(\Delta_{i} X)^2\leq r_h, (\Delta
_{i}M)^2\leq4r_h\}}\\
&&\qquad{}- \frac
{1}{\sqrt h} \sum
_{i=1}^n \biggl(\int_{t_{i-1}}^{t_i}\sigma_u \,\mathrm{d}W_u\biggr)h\, \mathrm{d}\bigl(2\sqrt[4]{
r_h}\bigr)I_{\{(\Delta_{i} X)^2\leq r_h, (\Delta_{i}M)^2
\leq4r_h\}},\nonumber
\end{eqnarray}
%
where, using the fact that $\int_{t_{i-1}}^{t_i}\sigma_u \,\mathrm{d}W_u$ and
$\Delta_{i} M^{(h)}$ are martingale
increments with zero quadratic covariation, the $L^1(\Omega)$-norm
of the first right-hand term is bounded by
\[
\sqrt{E\biggl[\frac{\sum_{i=1}^n (\int_{t_{i-1}}^{t_i}\sigma_u
\,\mathrm{d}W_u)^2(\Delta_{i} M^{(h)})^2}{h}\biggr]},
\]
%
which is dealt with similarly as in (\ref{sumDXzqDXdqsuh}) and tends to
zero. Moreover,
\begin{eqnarray*}
&&E\Biggl[\frac{1}{\sqrt h} \sum_{i=1}^n \biggl(\int
_{t_{i-1}}^{t_i}\sigma_u \,\mathrm{d}W_u
\biggr)h\, \mathrm{d}\bigl(2\sqrt[4]{ r_h}\bigr)I_{\{(\Delta_{i} X)^2\leq r_h, (\Delta
_{i}M)^2\leq4r_h\}}\Biggr]\\
&&\quad=
c\sqrt h \biggl[I_{\alpha\neq1}\bigl(c+r_h^{{(1-\alpha)}/{4}}\bigr)
+I_{\alpha=
1} \ln\frac{1}{r_h^{1/4}}\biggr]\\
&&\qquad {}\times E\Biggl[\sum_{i=1}^n \biggl(\int
_{t_{i-1}}^{t_i}\sigma_u
\,\mathrm{d}W_u\biggr)
I_{\{(\Delta_{i} X)^2\leq r_h, (\Delta_{i}M)^2\leq4r_h\}}\Biggr]\\
 &&\quad\leq
c\sqrt h \biggl[I_{\alpha\neq1}\bigl(c+r_h^{{(1-\alpha)}/{4}}\bigr)
+I_{\alpha=
1} \ln\frac{1}{r_h^{1/4}}\biggr]\sqrt{E\Biggl[\sum_{i=1}^n \biggl(\int
_{t_{i-1}}^{t_i}
\sigma_u \,\mathrm{d}W_u\biggr)^2\Biggr]} \rightarrow0.
\end{eqnarray*}
Using the fact that
\[
\frac{\sqrt{2h  \mathit{IQ}}}{\sqrt{2/3\sum
_i (\Delta_{i} X)^4
I_{\{ (\Delta_{i} X)^2\leq r_h\}} }}
\]
tends to 1 in probability,
treating $I_4(h)$ as in (\ref
{DoveDXqleqrhoDXdqleq4r}) and putting together the simplified version
of $I_2(h),$ we obtain that $(\hat{\mathit{IV}}_h-\mathit{IV})/\sqrt{2h\hat{\mathit{IQ}}_h}$
is the sum of a term which converges in distribution to an $N(0,1)$
r.v. plus a negligible term and a remainder
\begin{equation} \label{NormalizedBiasSimplified}
- \frac{\sum_{i=1}^n (\int_{t_{i-1}}^{t_i}\sigma_u \,\mathrm{d}W_u
)^2 I_{\{(\Delta_{i} X)^2>r_h
,(\Delta_{i} X_1)^2\leq4r_h\}} }{\sqrt{2h \mathit{IQ}}}
+\frac{ \sum_{i=1}^n (\Delta_{i}M)^2I_{\{(\Delta_{i} X)^2\leq r_h,
(\Delta_{i}M)^2\leq4r_h\}} }{\sqrt
{2h \mathit{IQ}}}.
\end{equation}

(a) If $\alpha<1$, the first term of (\ref
{NormalizedBiasSimplified}) is negligible
with respect to
$\frac{r_h^{1-\alpha/2}}{\sqrt{2h \mathit{IQ}}}$,
in fact,
\[
\frac{\sum_{i=1}^n (\int_{t_{i-1}}^{t_i}\sigma_u \,\mathrm{d}W_u)^2
I_{\{(\Delta_{i} X)^2>r_h,(\Delta_{i} X_1)^2
\leq4r_h\}}}{
r_h^{1-\alpha/2}}
\leq
\frac{\sum_{i=1}^n h\ln(1/ h )I_{\{(\Delta_{i} X)^2>r_h\}}}{
r_h^{1-\alpha/2}},
\]
where
\[
E\biggl[\frac{\sum_{i=1}^n h\ln(1/ h) I_{\{(\Delta_{i} X)^2>r_h\}}}{
r_h^{1-\alpha/2}}\biggr]\leq h^{1-\beta}\ln\frac1 h\rightarrow0.
\]
%
Therefore, (\ref{NormalizedBiasSimplified}) can be written as
\begin{equation} \label{NBulteriormSimlif}
\frac{r_h^{1-\alpha/2}}{\sqrt{2h \mathit{IQ}}}\biggl[\mathrm{o}_P(1) +
\frac{\sum_{i=1}^n (\Delta_{i}M)^2I_{\{(\Delta_{i} X)^2\leq r_h,
(\Delta_{i}M)^2\leq4r_h\}}}{r_h
^{1-\alpha/2}}\biggr].
\end{equation}
Using (\ref{I4}), Lemma \ref{StructDXdwhenleqr}(i) and Theorem \ref
{teoSpeedCvSumDXdq}, we arrive at
\begin{eqnarray*}
&&\frac{\sum_{i=1}^n (\Delta_{i}M)^2I_{\{(\Delta_{i} X)^2\leq r_h,
(\Delta_{i}M)^2\leq4r_h\}}}{r_h
^{1-\alpha/2}}\\
&&\quad\leq
\frac{\sum_{i=1}^n (\Delta_{i}M)^2I_{\{(\Delta_{i}M)^2\leq9r_h/4\}
}+ \mathrm{o}_P(h^{1-\alpha
/2})}{r_h^{1-\alpha/2}}\\
&&\quad\sim
\frac{\sum_{i=1}^n (\Delta_{i}M)^2I_{\{(\Delta_{i}M)^2\leq9r_h/4\}
}}{r_h^{1-\alpha/2}}\\
&&\qquad\leq
\frac{\sum_i (\int_{t_{i-1}}^{t_i}\!\int_{|x|\leq3\sqrt
{r_h}/2} x \tilde\mu
(\mathrm{d}x,\mathrm{d}t)-h\int_{3\sqrt{r_h}/2<|x|\leq1} x \nu(\mathrm{d}x) )^2
}{r_h^{1-\alpha
/2}}
\\
&&\quad= R_h + Tc + Tc\biggl(\frac h r_h\biggr)^{\alpha/2}h^{1-
\alpha/2} \stackrel{P}{\rightarrow} Tc,
\end{eqnarray*}
where the term $R_h$ has variance $\sim c r_h^{\alpha/2}\to0$ and so
converges to zero in probability. Since
$ \frac{r_h^{1-\alpha/2}}{\sqrt{h}} \rightarrow0,$ we arrive at
\[
\frac{\hat{\mathit{IV}}_h-\mathit{IV}}{\sqrt{2h\hat{\mathit{IQ}}_h}} \stackrel
{\mathit{st}}{\rightarrow} {N}(0,1).
\]

(b) If $\alpha> 1,$ define $R_t:= \sum_{s\leq t} I_{\{
|\Delta
M_{s}|>\sqrt h\}}$. Then, by (\ref{I4}), the last term (times $\sqrt{2
\mathit{IQ}}$) in~(\ref{NormalizedBiasSimplified}) dominates
\begin{eqnarray} \label{CLTalphamaguno}
&&\frac{\sum_{i=1}^n (\Delta_{i}M)^2I_{\{(\Delta_{i}M)^2\leq r_h/16\}
}- \mathrm{o}_P(h^{1-\alpha
/2})}{\sqrt{h}}
\nonumber\\
%
%
&&\quad=\frac{1}{\sqrt h} \biggl[\sum_i (\Delta_{i}M)^2I_{\{\Delta_i R=0\}}
+\sum_i (\Delta_{i}M)^2
\bigl[I_{\{(\Delta_{i}M)^2\leq r_h/16\}}- I_{\{\Delta_i R=0\}}\bigr]
\biggr] -
\mathrm{o}_P(h^{1/ 2-\alpha/2})
\nonumber
\\[-8pt]
\\[-8pt]
\nonumber
&&\quad\geq
- \mathrm{o}_P(h^{1 /2-\alpha/2})+
\frac{\sum_i  (\int_{t_{i-1}}^{t_i}\! \int_{|x|\leq\sqrt h}
x \tilde\mu
(\mathrm{d}x,\mathrm{d}t)-h\int_{\sqrt h<|x|\leq1} x \nu(\mathrm{d}x) )^2}{\sqrt h}\\
&&\qquad{}-
\frac
{\sum_i (\Delta_{i}M)^2I_{\{(\Delta_{i}M)^2>r_h/16, \Delta_i R=0\}
}}{\sqrt h}\nonumber
.
\end{eqnarray}
First,
\begin{eqnarray*}
&&\sum_i (\Delta_{i}M)^2I_{\{(\Delta_{i}M)^2>{r_h}/{16}, \Delta
_i R=0\}}\\
&&\quad= \sum_i
\biggl[\Delta_i [M] +
2\int_{t_{i-1}}^{t_i}(M_{s_{-}}- M_{t_{i-1}}) \,\mathrm{d}M_{s}\biggr] I_{\{
(\Delta_{i}M)^2>{r_h
}/{16}, \Delta_i R=0\}}.
\end{eqnarray*}
As in Lemma \ref{StructDXdwhenleqr}, the sum of the right-hand terms
within brackets is of order $u_n= (n/\log n)^{1/\alpha}$ so that
\[
\frac{\sum_i|\int_{t_{i-1}}^{t_i}(M_{s_{-}}- M_{t_{i-1}})
\,\mathrm{d}M_{s}|}{\sqrt h}=
\frac{u_n\sum_i|\int_{t_{i-1}}^{t_i}(M_{s_{-}}- M_{t_{i-1}})
\,\mathrm{d}M_{s}|}{u_n\sqrt
h}\stackrel{P}{\rightarrow} 0
\]
since $u_n\sqrt h= (\frac{n^{(1-\alpha/2)}}{\log n}
)^{1/ \alpha}\rightarrow+\infty$.
Theorem \ref{teoSpeedCvSumDXdq} applied with $u=1/2$ yields that with
$\varepsilon=h^{1/ 2},$
\[
\sum_i\biggl(\int_{t_{i-1}}^{t_i}\!\int_{|x|\leq\sqrt h} x \tilde\mu
(\mathrm{d}x,\mathrm{d}t)-h\int
_{\sqrt h<|x|\leq1} x \nu(\mathrm{d}x) \biggr)^2= \varepsilon^{2-
\alpha/2}Y_h +
Tc\varepsilon^{2-\alpha} + Tc h \varepsilon^{2 - 2\alpha},
\]
where $\operatorname{var}(Y_h)\to1$. Therefore, in (\ref{CLTalphamaguno}), we remain
with
\[
h^{1 /2 - \alpha/2} \biggl[- \mathrm{o}_P(1) +h^{\alpha/4}Y + Tc
+ Tc h^{1-\alpha/2}
-\frac{\sum_i \Delta_i [M] I_{\{(\Delta_{i}M)^2>r_h/16, \Delta_i
R=0\}} }{h^{1
-\alpha/2}}\biggr]\stackrel{\mathit{a.s.}}{\rightarrow} +\infty,
\]
%
where the divergence is due to the fact that $h^{1/ 2 -
\alpha/2}\rightarrow+\infty$ while
$\frac{\sum_i \Delta_i [M] I_{\{(\Delta_{i}M)^2>r_h/16, \Delta_i
R=0\}} }{h^{1
-\alpha/2}}$ tends to zero in probability since its expected
value is dominated by
\begin{eqnarray*}
&&\frac{n}{h^{1 -\alpha/2}} E^{1 /2}\bigl[\bigl(\Delta_i [M]
I_{\{ \Delta_i R=0\}} \bigr)^2\bigr] P^{1/ 2}\{(\Delta
_{i}M)^2>r_h/16, \Delta
_i R=0\}
\\
&&\quad\leq\frac{n}{h^{1 -\alpha/2}} \biggl(h \int_{|x|\leq\sqrt h }
x^{4} \nu(\mathrm{d}x)\biggr)^{1/ 2} h^{(2-\alpha/2 - \beta)1/2 }
= h^{ {(1-\beta)}/{2}}\rightarrow0,
\end{eqnarray*}
having used the fact that
\begin{eqnarray} \label{PDXdgrmasaltisottosoglia}P\{(\Delta
_{i}M)^2>r_h, \Delta_i R=0\}&=&P\bigl\{
(\Delta_{i}M)^2I_{\{ \Delta_i R=0\}}>r_h\bigr\} \leq\frac{E[(\Delta
_{i}M)^2I_{\{ \Delta_i
R=0\}}]}{r_h}
\nonumber
\\[-8pt]
\\[-8pt]
\nonumber
&=& \frac{h\int_{|x|\leq\sqrt h} x^2 \nu(\mathrm{d}x)}{r_h} = h^{2-
\alpha/2
-\beta}.
\end{eqnarray}
%
On the other hand, the first term in (\ref{NormalizedBiasSimplified})
is negligible with respect to $h^{1/ 2-\alpha/2}$ (the speed
of divergence of $(\sum_{i=1}^n (\Delta_{i}M)^2I_{\{(\Delta
_{i}M)^2\leq r_h/16\}}-
\mathrm{o}_P(h^{1-\alpha/2}))/\sqrt{h}$) because
\[
\frac{\sum_{i=1}^n (\int_{t_{i-1}}^{t_i}\sigma_u \,\mathrm{d}W_u)^2
I_{\{(\Delta_{i} X)^2>r_h,(\Delta_{i} X_1)^2
\leq4r_h\}}
}{\sqrt h h^{1/ 2 - \alpha/2}}\leq\frac{h \log(1/ h)
h^{-{\alpha\beta}/{2}}
}{h^{1 -\alpha/2}}= h^{(\alpha/2) (1-\beta)}\log\frac1 h
\rightarrow0.
\]
%
Therefore, (\ref{NormalizedBiasSimplified}) explodes to $+\infty.$
Finally, if $\alpha=1$ in (\ref{NormalizedBiasSimplified}), then the first
term is negligible,~as
\[
\frac{\sum_{i=1}^n (\int_{t_{i-1}}^{t_i}\sigma_u \,\mathrm{d}W_u)^2
I_{\{(\Delta_{i} X)^2>r_h,(\Delta_{i} X_1)^2
\leq4r_h\}}
}{\sqrt h}= \mathrm{O}_p\biggl(h^{{(1-\beta)}/{2}} \log\frac1 h\biggr) \rightarrow0.
\]
For the second term, we take a $\delta>0$ such that $2/3<\beta+\delta
<1$, we choose $\varepsilon=h^{{(\beta+\delta)}/{2}}$ and we use
the same
steps as were used to reach (\ref{CLTalphamaguno}) for $\alpha>1$, but
we consider 
$\tilde R_t=\sum_{s\leq t} I_{\{|\Delta M_{s}|>\varepsilon\}}$ in
place of
$R_t$. Also using Theorem \ref{teoSpeedCvSumDXdq}, we obtain that the
second term in \eqref{NormalizedBiasSimplified} dominates
\begin{eqnarray*}
&&\frac{Y_h \varepsilon^{3/ 2}}{\sqrt{2h \mathit{IQ}}} +\frac{\varepsilon
}{\sqrt{2h \mathit{IQ}}}
- \frac{\sum\Delta_i[M] I_{\{(\Delta_{i}M)^2>{r_h}/{16},
\Delta_i \tilde
R=0\}}}{\sqrt{2h \mathit{IQ}}}\\
&&\quad{}-
\frac{2 \sum_{i=1}^n \int_{t_{i-1}}^{t_i}(M_{ s_{-}} - M_{ t_{i-1}})
\,\mathrm{d}M_{ s}\,
I_{\{(\Delta_{i}M)^2>{r_h}/{16}, \Delta_i \tilde R=0\}} }{\sqrt{2h
\mathit{IQ}}},
\end{eqnarray*}
where the variance of $Y_h$ tends to 1 so that $Y_h
\varepsilon^{3/2}/\sqrt{h}$ tends to zero in probability. The second term tends to
$+\infty$ at rate $\varepsilon/\sqrt h$.
The third term is negligible with respect to $\varepsilon/\sqrt h$: applying
(\ref{PDXdgrmasaltisottosoglia}) with $\tilde R$ in place of $R$ and
the Cauchy--Schwarz inequality, we get
\[
E\biggl[\frac{1}{\varepsilon}\sum\int_{t_{i-1}}^{t_i}\!\int_{|x|\leq
1} x^2 \mu(\mathrm{d}x, \mathrm{d}t) I_{\{
(\Delta_{i}M)^2>{r_h}/{16}, \Delta_i \tilde R=0\}}\biggr]=
\mathrm{O}(h^{{\delta}/{2}})\rightarrow0.
\]
Finally, the last term is also negligible since the speed of
convergence to zero of the numerator is $u_n= n/\log^2 n$ (as in the
proof of Lemma \ref{StructDXdwhenleqr}) and
$u_n\sqrt h\rightarrow+\infty$.
So, even for $\alpha=1,$ the normalized bias $(\hat{\mathit{IV}}_h -
\mathit{IV})/\sqrt{2h\hat{\mathit{IQ}}_h}$ diverges to $+\infty$.
\end{pf*}

\begin{pf*}{Proof of Proposition \protect\ref{EstimJwhenAlphaLess1}} As in Lemma
\ref{IndDxqleqrvsIndDxdqleqr} with $\sqrt{r_h}$ in place of $r_h$ as
bound for $\max_{i=1,\ldots,n} |a_{ni}|$, using the fact that $\alpha
<1$ and
applying Lemma \ref{StructDXdwhenleqr}(i), we deduce that $\hat
H_{h}$ has the same limit in probability as
\[
X_T-\sum_{i=1}^n (\Delta_{i} X_0+ \Delta_{i}M) I_{\{\Delta_{i} N=0,
(\Delta_{i}M)^2\leq r_h\}}
\]
when $h\rightarrow0$. Moreover, since a.s. $N_T<\infty$ and $\sum_{i=1}^n
\Delta_{i} X_0I_{\{(\Delta_{i}M)^2> r_h\}})=\mathrm{O}_P(h^{(1-\alpha\beta
)/2}\times \sqrt{\log
(1/h)})\rightarrow0$, taking $\tilde R_t = \sum_{s\leq t} I_{\{
|\Delta
M_{s}|>\sqrt{r_h}\}},$
the above term has limit in probability equal to
\begin{eqnarray*}
&&X_T-\lim_h \sum_{i=1}^n \bigl(\Delta_{i} X_0+ \Delta_{i}MI_{\{(\Delta
_{i}M)^2\leq r_h\}}\bigr)\\
&&\quad=
X_T-X_{0 T} -
\lim_h \Biggl[\sum_{i=1}^n \int_{t_{i-1}}^{t_i}\!\int_{|x|\leq\sqrt
{r_h}} x \tilde\mu(\mathrm{d}x,\mathrm{d}t)
-T\int_{\sqrt{r_h}<|x|\leq1} x \nu(\mathrm{d}x) \Biggr]\\
&&\qquad{}-\lim_h \sum_i\Delta
_{i}M\bigl(I_{\{(\Delta_{i}M)^2
\leq r_h\}}- I_{\{\Delta_i \tilde R =0\}}\bigr).
\end{eqnarray*}
Using the fact that $P\{\Delta_i \tilde R\geq1\}= \mathrm{O}(h^{1-\alpha\beta
/2}),$ as was used after (\ref{ProbDtildeNgeq1}), we deduce that $\sum
_i\Delta_{i}MI_{\{(\Delta_{i}M)^2\leq r_h, \Delta_i \tilde R \geq1\}
}=\mathrm{O}_P(h^{(1-\alpha
)\beta/2})$ $\rightarrow0$. Using the H\"older inequality with exponents
$p=q=2$, we have
$\sum_i\Delta_{i}MI_{\{(\Delta_{i}M)^2> r_h, \Delta_i \tilde R =0\}
}=\mathrm{O}_P(r_h
^{(1-\alpha)\beta/2})\rightarrow0$. Finally,\vspace*{-2pt} $\int_0^T\!\!\!\int
_{|x|\leq\sqrt{r_h}} x
\tilde\mu(\mathrm{d}x, \mathrm{d}t)
\stackrel{L^2}{\rightarrow} 0$ and $\int_{\sqrt{r_h}< |x|\leq1} x
\nu(\mathrm{d}x)\rightarrow m $
so that $\hat H_{h, T}\stackrel{P}{\rightarrow} J_T +mT$.
\end{pf*}
\end{appendix}

\section*{Acknowledgements} A previous version of this working paper
appeared as ``Nonparametric test for analyzing the fine structure of
price fluctuations'', Columbia Financial Engineering Report 2007-13.
This research was supported in part by the European Science Foundation
program ``Advanced Mathematical Methods in Finance'', by Istituto
Nazionale di Alta Matematica and by MIUR Grants Nos 206132713-001 and
2004011204-002. We thank Jean Jacod and Suzanne Lee for important comments.

\printhistory

\end{document}